\documentclass[12pt]{amsart}
\usepackage[all]{xy}
\usepackage{a4wide}
\usepackage{mathtools}
\usepackage{amssymb}
\usepackage{amsthm}
\usepackage{amsmath}
\usepackage{amscd}
\usepackage{verbatim}
\usepackage[all]{xy}
\usepackage{hyperref}

\usepackage{comment}

\usepackage[obeyDraft,textsize=tiny]{todonotes}

\usepackage{enumitem}

\numberwithin{equation}{section}

\DeclareMathOperator{\GL}{GL}
\DeclareMathOperator{\id}{id}
\DeclareMathOperator{\G}{\Gamma}

\DeclareMathOperator{\SL}{SL}
\DeclareMathOperator{\Mp}{Mp}
\DeclareMathOperator{\Sp}{Sp}

\DeclareMathOperator{\N}{\mathbb{N}}
\DeclareMathOperator{\Z}{\mathbb{Z}}
\DeclareMathOperator{\R}{\mathbb{R}}
\DeclareMathOperator{\C}{\mathbb{C}}
\DeclareMathOperator{\Q}{\mathbb{Q}}

\renewcommand{\H}{\mathbb{H}}

\DeclareMathOperator{\tr}{tr}

\DeclareMathOperator{\cusp}{cusp}
\DeclareMathOperator{\mero}{mero}

\DeclareMathOperator{\e}{\mathfrak{e}}
\newcommand{\Gr}{\mathrm{Gr}}
\newcommand{\reg}{\mathrm{reg}}

\newcommand{\Iso}{\mathrm{Iso}}
\newcommand{\CT}{\mathrm{CT}}
\newcommand{\SO}{\mathrm{SO}}
\newcommand{\Orth}{\mathrm{O}}

\newcommand{\bs}{\ensuremath{\backslash}}

\newcommand{\omegamero}[1]{\omega_{#1}^\text{mero}}

\newcommand{\omegacusp}[1]{\omega_{#1}^\text{cusp}}

\newcommand{\avgOmegacusp}[2]{\Omega^\text{cusp}_{#1, #2}}
\newcommand{\avgomegacusp}[2]{\omega_{#1, #2}^\text{cusp}}

\newcommand{\DN}{\mathrm{DN}}
\newcommand{\Millson}{\mathrm{M}}
\newcommand{\Shintani}{\mathrm{Sh}}
\newcommand{\Oda}{\mathrm{Oda}}
\newcommand{\Pet}{\mathrm{Pet}}

\makeatletter
\newcommand*\defbb[1]{
  \expandafter\newcommand\csname I#1\endcsname{\mathbb{#1}}}
\newcommand*\defbbs[1]{
  \@for\@i:=#1\do{\expandafter\defbb\expandafter{\@i}}}
\makeatother
\defbbs{A,B,C,D,E,F,G,H,I,K,L,M,N,O,P,Q,R,S,T,U,V,W,X,Y,Z} 

\makeatletter
\newcommand*\deffrak[1]{
  \expandafter\newcommand\csname frak#1\endcsname{\mathfrak{#1}}}
\newcommand*\deffraks[1]{
  \@for\@i:=#1\do{\expandafter\deffrak\expandafter{\@i}}}
\makeatother
\deffraks{a,b,c,d,e,f,g,h,i,j,k,l,m,n,o,p,q,r,s,t,u,v,w,x,y,z}

\makeatletter
\newcommand*\defcal[1]{
  \expandafter\newcommand\csname cal#1\endcsname{\mathcal{#1}}}
\newcommand*\defcals[1]{
  \@for\@i:=#1\do{\expandafter\defcal\expandafter{\@i}}}
\makeatother
\defcals{A,B,C,D,E,F,G,H,I,J,K,L,M,N,O,P,Q,R,S,T,U,V,W,X,Y,Z}

  \newtheorem{Theorem}{Theorem}[section]
  \newtheorem{Lemma}[Theorem]{Lemma}
  \newtheorem{Proposition}[Theorem]{Proposition} 
  \newtheorem{Corollary}[Theorem]{Corollary}
  
  \theoremstyle{definition} 
  \newtheorem{Definition}[Theorem]{Definition}
  
  \newtheorem{Remark}[Theorem]{Remark}

\title{Cycle integrals of meromorphic Hilbert modular forms}

\author{Claudia Alfes}

\address{Universität Bielefeld, Fakultät für Mathematik, Postfach 100 131, 33501 Bielefeld, Germany}
\email{alfes@math.uni-bielefeld.de}

\author{Baptiste Depouilly}

\address{ETH, Mathematics Dept., CH-8092, Z\"urich, Switzerland}
\email{baptiste.depouilly@math.ethz.ch}

\author{Paul Kiefer}

\address{Universität Bielefeld, Fakultät für Mathematik, Postfach 100 131, 33501 Bielefeld, Germany}
\email{pkiefer@math.uni-bielefeld.de}

\author{Markus Schwagenscheidt}

\address{ETH, Mathematics Dept., CH-8092, Z\"urich, Switzerland}

\email{mschwagen@ethz.ch}

\thanks{The research of the C. Alfes was supported by the Daimler and Benz Foundation.}
\thanks{The research of the C. Alfes and P. Kiefer is funded by the Deutsche Forschungsgemeinschaft (DFG, German Research Foundation) -- SFB-TRR 358/1 2023 -- 491392403}
\thanks{The research of M. Schwagenscheidt is supported by SNF grant PZ00P2\_202210}

\begin{document}

\begin{abstract}
We establish a rationality result for linear combinations of traces of cycle integrals of certain meromorphic Hilbert modular forms. These are meromorphic counterparts to the Hilbert cusp forms $\omega_m(z_1,z_2)$, which Zagier investigated in the context of the Doi-Naganuma lift. 
We give an explicit formula for these cycle integrals, expressed in terms of the Fourier coefficients of harmonic Maass forms.
A key element in our proof is the explicit construction of locally harmonic Hilbert-Maass forms on $\H^2$, which are analogous to the elliptic locally harmonic Maass forms examined by Bringmann, Kane, and Kohnen. Additionally, we introduce a regularized theta lift that maps elliptic harmonic Maass forms to locally harmonic Hilbert-Maass forms and is closely related to the Doi-Naganuma lift.
\end{abstract}

\maketitle

\setcounter{tocdepth}{1} 
\tableofcontents

\section{Introduction and statement of results}
\label{sec:introduction}

In \cite{zagierdoinaganuma} Zagier constructed Hilbert cusp forms for the group $\SL_2(\mathcal{O}_F)$, where $\mathcal{O}_F$ is the ring of integers of a real quadratic field $F$. They are defined for natural numbers $m > 0$ and even integers $k \geq 4$ by
\begin{align}\label{omegamcusp}
\omega_m(z_1,z_2) = \sum_{\substack{a,b \in \Z, \, \nu \in \mathfrak{d}_F^{-1} \\ N(\nu) - ab = m/D}}\frac{1}{(az_1 z_2 + \nu z_1 + \nu' z_2 + b)^{k}},
\end{align}
where $D > 0$ is the discriminant of $F$, $\mathfrak{d}_F = (\sqrt{D})$ is the different, $N(\nu) = \nu\nu'$ is the norm and $\nu'$ the conjugate of $\nu$. These functions define Hilbert cusp forms of parallel weight $k$ for $\SL_2(\mathcal{O}_F)$. The main result of \cite{zagierdoinaganuma} states that the generating series of the cusp forms $\omega_m(z_1,z_2)$ is a kernel function for the Doi-Naganuma correspondence between elliptic modular forms of weight $k$ for $\Gamma_0(D)$ with Nebentypus $\chi_D = \left( \frac{D}{\cdot}\right)$ and Hilbert modular forms of weight $k$ for $\SL_2(\mathcal{O}_F)$.

If we allow $m$ to be negative in the definition \eqref{omegamcusp}, we obtain \emph{meromorphic} Hilbert modular forms of weight $k$. In order to clearly distinguish between the two cases, we will denote the cusp forms by $\omega_m^{\cusp}(z_1,z_2)$ for $m > 0$ and the meromorphic Hilbert modular forms by $\omega_n^{\mero}(z_1,z_2)$ for $n < 0$. The main goal of the current work is to investigate the cycle integrals of these meromorphic modular forms $\omega_n^{\mero}(z_1,z_2)$ along certain real analytic cycles in $\H^2$, and in particular to study their rationality properties. Our work is motivated by analogous rationality results for elliptic modular forms. Namely, the Hilbert modular forms $\omega_m^{\cusp}$ and $\omega_n^{\mero}$ can be viewed as the natural analogues of the elliptic modular forms 
\[
f_{k,d}(z) = \sum_{\substack{a,b,c \in \Z \\ b^2-4ac = d}}\frac{1}{(az^2 + bz + c)^{k}},
\]
which have also been constructed for the first time in \cite{zagierdoinaganuma}. Indeed, the generating series of the cusp forms $f_{k,d}$ for $d > 0$ is a kernel function for the Shimura-Shintani correspondence \cite{kohnenzagierwaldspurger} between elliptic modular forms of integral and half-integral weight. The rationality of the cycle integrals of the cusp forms $f_{k,d}$ for $d > 0$ is a classical result due to Kohnen and Zagier \cite{kohnenzagier} and Katok \cite{katok}, while the analogous rationality questions for the cycle integrals of the meromorphic modular forms $f_{k,d}$ for $d < 0$ have been studied more recently, see \cite{anbs, loebrichschwa2, anbms, loebrichschwa1}.

We now describe our results in more detail.
We first introduce some notation. Consider the lattice
\[
L = \left\{\begin{pmatrix}a & \nu' \\ \nu & b \end{pmatrix}\, : \, a,b \in \Z, \,  \nu \in \mathcal{O}_F \right\}
\]
together with the quadratic form $Q(X) = -\det(X)$. Its dual lattice $L'$ is given by matrices of the same shape, with $a,b \in \Z$, but $\nu \in \mathfrak{d}_F^{-1}$. For $m \in \Z$ we let $L_m$ be the set of all $X \in L'$ with $Q(X) = m/D$. The group $\Gamma = \SL_2(\mathcal{O}_F)$ acts isometrically on $L$ by $\gamma.X = \gamma X \gamma'^t$. For $m \neq 0$ it acts on $L_m$ with finitely many orbits. Throughout, we will write $Z = (z_1,z_2) \in \H^2$ for brevity. For $X = \left(\begin{smallmatrix}a & \nu' \\ \nu & b \end{smallmatrix} \right) \in L'$ we define the quantities
\begin{align*}
q_Z(X) &= -bz_1 z_2 + \nu z_1 + \nu' z_2 - a, \\
p_Z(X) &= \frac{1}{y_1}(-b\overline{z}_1 z_2 + \nu \overline{z}_1  + \nu' z_2 -a).
\end{align*}
Using this notation, the meromorphic Hilbert modular forms for $n < 0$ can be written as
\[
\omega_n^{\mero}(Z) = \sum_{X \in L_n}q_Z(X)^{-k}.
\]
They have poles along the \emph{algebraic cycles}
\[
T_n = \bigcup_{X \in L_n}T_X, \qquad T_X = \{Z \in \H^2 \, :\, q_Z(X) = 0\},
\]
which are also called \emph{Hirzebruch-Zagier divisors} of discriminant $n$.

For a matrix $Y \in L'$ with $Q(Y) > 0$ we consider the \emph{real analytic cycle}
\[
C_Y = \{Z \in \H^2: p_Z(Y) = 0\},
\]
which is real two-dimensional. For a Hilbert cusp form $F(Z)$ of weight $k$ for $\Gamma$ and $m > 0$ we define its \emph{$m$-th trace of cycle integrals} as the finite sum
\[
\tr_m(F) = \sum_{Y \in \Gamma \backslash L_m}\int_{\Gamma_Y \backslash C_Y}F(Z) \, q_Z(Y)^{k-2} dz_1 dz_2,
\]
where $\Gamma_Y$ denotes the stabilizer of $Y$ in $\Gamma$. Note that $q_Z(Y)$ transforms with parallel weight $-1$ under $\Gamma_Y$, so the cycle integral is well-defined. 

Now, for $n < 0$ the meromorphic Hilbert modular form $\omega_n^{\mero}(Z)$ decays like a cusp form, but it has poles along the algebraic cycle $T_n$. For simplicity, let us assume that $T_n$ does not intersect any of the real analytic cycles $C_Y$ for $Y \in L_m$. Then the trace of cycle integrals $\tr_m(\omega_n^{\mero})$ is well-defined. Moreover, in order to simplify the statement of our main result, we assume in the introduction that the discriminant $D$ of $F$ is an odd prime, and that a certain space of cusp forms vanishes. Namely, we consider the space $S_{k}^+(\Gamma_0(D),\chi_D)$ of cusp forms of weight $k$ and character $\chi_D$ for $\Gamma_0(D)$ whose Fourier coefficients satisfy the ``plus space'' condition $a(n) = 0$ if $\chi_D(n) = -1$. Then our main result is the following rationality statement.

\begin{Theorem}\label{thm:MainResultIntro}
	Let $k \geq 4$ be an even integer and let $m> 0$ and $n < 0$. Suppose that the discriminant $D$ of the real quadratic field $F$ is an odd prime, and that the space $S_{k}^+(\Gamma_0(D),\chi_D)$ is trivial. Moreover, assume that $T_n$ does not intersect any of the cycles $C_Y$ for $Y \in L_m$. Then the traces of cycle integrals 
	\[
	\tr_m(\omega_n^{\mero}) = \sum_{Y \in \Gamma \backslash L_m}\int_{\Gamma_Y \backslash C_Y}\omega_n^{\mero}(Z) \, q_Z(Y)^{k-2} dz_1 dz_2
	\]
	of the meromorphic Hilbert modular forms $\omega_n^{\mero}$  are \emph{rational} multiples of $\pi i$.
\end{Theorem}

This result is an analogue of \cite[Theorem~1.1]{anbs}, where the rationality of the traces of cycle integrals of the elliptic meromorphic modular forms $f_{k,d}$ for $d < 0$ was proved. In the body of the paper we will allow arbitrary fundamental discriminants $D$ by working with vector-valued modular forms for the Weil representation associated with an even lattice of signature $(2,2)$. Moreover, we will remove the restriction on $S_{k}^+(\Gamma_0(D),\chi_D)$ by taking suitable linear combinations of traces of cycle integrals. We state our main result in the general setting in Theorem~\ref{thm:MainRationalityResult}.

The rationality statement of Theorem~\ref{thm:MainResultIntro} follows from an explicit formula for $\tr_m(\omega_n^{\mero})$ in terms of Fourier coefficients of harmonic Maass forms and unary theta functions, see Theorem~\ref{thm:ExplicitFormula}. We will give a sketch of the proof at the end of the introduction.  A major step in the proof is the construction of \emph{locally harmonic Hilbert-Maass forms}, which are natural analogues of the elliptic locally harmonic Maass forms constructed by H\"ovel \cite{hoevel} and Bringmann, Kane, and Kohnen \cite{brikakohnen}. However, in the Hilbert modular surface case it is more natural to define them as differential forms on $\Gamma \backslash\H^2$.

\begin{Definition}\label{def:LocallyHarmonicMaassForms}
For a natural number $m > 0$ and an even integer $k \geq 4$ we define the \emph{locally harmonic Hilbert-Maass form}
\[
\Omega_m^{\cusp}(Z) = \Omega_{m,1}^{\cusp}(Z)\, dz_1 \wedge d\mu(z_2) + \Omega_{m,2}^{\cusp}(Z)\, d\mu(z_1) \wedge dz_2,
\]
with the usual invariant measure $d\mu(z) = \frac{dx \, dy}{y^2}$ on $\H$, and the functions
\[
\Omega_{m,1}^{\cusp}(Z) = \sum_{Y \in L_m}y_1^{k-2}y_2^k\frac{\overline{q_Z(Y)}^{1-k}}{\overline{p_Z(Y)}}, \qquad \Omega_{m,2}^{\cusp}(Z) = \sum_{Y \in L_m}y_1^{k}y_2^{k-2}\frac{\overline{q_Z(Y)}^{1-k}}{\frac{y_1}{y_2}p_Z(Y)}.
\]
\end{Definition}

The functions $\Omega_{m,1}^{\cusp}(Z)$ and $\Omega_{m,2}^{\cusp}(Z)$ transform like Hilbert modular forms of weight $(2-k,-k)$ and $(-k,2-k)$ for $\Gamma$, respectively. This implies that $\Omega_{m}^{\cusp}(Z)$ defines a $(2,1)$-form on $\Gamma \backslash \H^2$ with values in the sheaf $\mathcal{L}_{-k}$ of functions on $\H^2$ that transform like Hilbert modular forms of parallel weight $-k$. Note that the elliptic locally harmonic Maass forms constructed in \cite{brikakohnen} have jumps along certain closed geodesics in $\H$, whereas our locally harmonic Hilbert-Maass forms blow up along the real analytic cycles $C_Y$ for $Y \in L_m$. 

The functions $\Omega_{m,1}^{\cusp}(Z)$ and $\Omega_{m,2}^{\cusp}(Z)$ are related to the cusp form $\omega_m^{\cusp}(Z)$ via the $\xi$-operator $\xi_{k} = 2iy^{k} \overline{\frac{\partial}{\partial \overline{z}}}$ from the theory of harmonic Maass forms. A straightforward calculation shows that we have
\begin{align}
\begin{split}\label{xiequation}
\xi_{2-k,z_1}\Omega_{m,1}^{\cusp}(Z) &= (k-1)y_2^k\, \omega_m^{\cusp}(Z), \\
\xi_{2-k,z_2}\Omega_{m,2}^{\cusp}(Z) &= (k-1)y_1^k\, \omega_m^{\cusp}(Z).
\end{split}
\end{align}
In Section~\ref{sec:XiOperator} we define an analogue of the $\xi$-operator in the Hilbert case. It acts on differential forms on Hilbert modular surfaces. The relations \eqref{xiequation} then imply the following result.

\begin{Proposition}\label{prop:XiRelationIntro}
	We have
	\[
	\xi_{-k,Z}\Omega_{m}^{\cusp}(Z) = -2(k-1) \, \omega_m^{\cusp}(Z),
	\]
	where $\xi_{-k,Z}$ is the differential operator on $(2,1)$-forms defined in Section~\ref{sec:XiOperator}. Moreover, $\Omega_{m}^{\cusp}(Z)$ is harmonic outside of its singularities.
\end{Proposition}

\begin{Remark}
	We define \emph{polar harmonic Hilbert-Maass forms} $\Omega_{n}^{\mero}(Z)$ for $n < 0$ by the same formula as in Definition~\ref{def:LocallyHarmonicMaassForms}. They are the analogues of the elliptic polar harmonic Maass forms studied by Bringmann and Kane \cite{bringmannkanepolar}, and they are related to $\omega_n^{\mero}$ by the $\xi$-operator by the same formula as in Proposition~\ref{prop:XiRelationIntro}. 
\end{Remark}

 In Section~\ref{sec:ThetaLifts} we will show that $\Omega_m^{\cusp}(Z)$ can be written as a regularized theta lift of a harmonic Maass form. To this end, we define the \emph{Millson theta lift} of a vector-valued harmonic Maass form $g(\tau)$ of weight $2-k$ for the dual Weil representation of $L$ by
 \begin{align}\label{theta lift representation}
 \Phi_L^{\Millson}(g,Z) =  \int_{\SL_2(\Z) \backslash \H}^{\reg} \left\langle g(\tau),\overline{\Theta_L^{\Millson}(\tau,Z)} \right\rangle v^{2-k} d\mu(\tau), 
 \end{align}
 where the integral has to be regularized as explained by Borcherds \cite{borcherds}, $\langle \cdot,\cdot \rangle$ denotes the natural inner product on the group ring $\C[L'/L]$, and $\Theta_L^{\Millson}(\tau,Z)$ is the \emph{Millson theta form} as defined in Section~\ref{sec:O(2,2)theta}. The Millson theta form is a $(2,1)$-form on $\Gamma \backslash \H^2$ with values in $\mathcal{L}_{-k}$ as a function of $Z$, and its complex conjugate transforms like a vector-valued modular form of weight $2-k$ for the dual Weil representation of $L$ in $\tau$. We let $f_{m}(\tau)$ be the vector-valued Maass-Poincar\'e series of weight $2-k$ and index $m$ as defined in Section~\ref{subsec:MaassPoincareSeries}. By the unfolding argument, we obtain the following result.
 
 \begin{Theorem}[Theorem~\ref{thm:unfoldpc}]\label{thm:ThetaLiftIntro}
 	We have
	\[
	\Phi_L^{\Millson}(f_{m},Z) = \frac{2}{\pi} (4m/D)^{k - 1}\,\Omega_m^{\cusp}(Z).
	\]
 \end{Theorem}

 The above theta lift is related to the Doi-Naganuma theta lift $\Phi_L^{\DN}(f,Z)$ of a cusp form $f$ of weight $k$ (defined in Section~\ref{sec:ThetaLifts}) via the $\xi$-operator. 
  \begin{Theorem}[Theorem~\ref{thm:ThetaLiftAndXiOperator}]\label{xidiagram}
  For a harmonic Maass form $g$ of weight $2-k$ for the dual Weil representation of $L$ we have
  \[
  \xi_{-k,Z}\Phi_L^{\Millson}(g,Z) = -2\Phi_L^{\DN}(\xi_{2-k,\tau}g,Z).
  \]
  \end{Theorem}
  
  Using this relation, in Theorem~\ref{thm:DistributionEquation} we show that the Millson lift satisfies a current equation analogously to the current equation for the Millson lift in signature $(2,1)$, which has been studied by H\"ovel \cite{hoevel} and Crawford and Funke \cite{crawfordfunke}.

We utilize the locally harmonic and polar harmonic Hilbert-Maass forms $\Omega_{m}^{\cusp}$ and $\Omega_n^{\mero}$ to evaluate regularized Petersson inner products of $\omega_{m}^{\cusp}$ and $\omega_{n}^{\mero}$ against cusp forms. In fact, we derive a general result on the evaluation of regularized Petersson inner products
  \[
  \langle G,\xi_{-k,Z}H \rangle_{\Pet}^{\reg} = \int_{\Gamma \bs \IH^2}^{\reg} F(Z) \overline{\xi_{-k,Z} H(Z)} (y_1 y_2)^k d\mu(Z)
  \]
  of Hilbert modular forms $F$ and $G = \xi_{-k,Z}H$ of weight $k$, where $F$ has singularities along algebraic cycles and $H$ has singularities along real cycles in Theorem~\ref{thm:innerProductEvaluation} of specific types (compare~\ref{subsec:RegularizedInnerProducts}). Here $d\mu(Z) = d\mu(z_1)\, d\mu(z_2)$, and the inner product has to be regularized as explained in Definition~\ref{def:RegularizedInnerProduct}.

  \begin{Theorem}\label{thm:InnerProductFormulas}
  	Let $F$ be a Hilbert cusp form of even weight $k \geq 4$ for $\Gamma$. Let $m > 0$ and $n < 0$. Then we have
  	\[
	\langle F,\omega_m^{\cusp}\rangle_{\Pet}^{\reg} = \frac{2\pi i}{(4m/D)^{k-1}(k-1)}\tr_m(F) \qquad \text{and} \qquad \langle F, \omega_n^{\mero} \rangle_{\Pet}^{\reg} = 0.
	\]
  \end{Theorem}
  
  The first formula is well-known and has been proved in \cite{ZagierCycleIntegralsOmega}, Theorem~6. We will give a generalization and a different proof in Theorem~\ref{thm:innerProductEvaluation}. The second formula involving $\omega_n^{\mero}$ is proved in Theorem~\ref{thm:orthogonalToCuspForms}.

We are now in the position to sketch the proof of Theorem~\ref{thm:MainResultIntro}:
\begin{enumerate}
	\item We have seen in Theorem~\ref{thm:InnerProductFormulas} that $\omega_n^{\mero}$ is orthogonal to cusp forms. In particular, it is orthogonal to $\omega_m^{\cusp}$. On the other hand, we can evaluate the regularized Petersson inner product $\langle \omega_m^{\cusp},\omega_n^{\mero}\rangle_{\Pet}^{\reg}$ using Proposition~\ref{prop:XiRelationIntro} and Stokes' Theorem. The result is essentially given by $\tr_m(\omega_n^{\mero})+\tr_{n}(\Omega_{m}^{\cusp})$, where the second trace is a sum of ``evaluations'' of $\Omega_{m}^{\cusp}$ at the algebraic cycles $T_X$ over $X \in \Gamma \backslash L_n$ (see Theorem~\ref{thm:innerProductEvaluation}). Summarizing, we obtain the ``duality'' 
	\[
	\tr_m(\omega_n^{\mero}) \stackrel{\cdot}{=} \tr_n(\Omega_{m}^{\cusp}),
	\]
	where $\stackrel{\cdot}{=}$ means equality up to a non-zero constant multiple.
	\item In order to evaluate $\tr_n(\Omega_m^{\cusp})$, we write $\Omega_m^{\cusp}$ as a Millson theta lift of a harmonic Maass form $f_{m}$ as in Theorem~\ref{thm:ThetaLiftIntro}. This leads us to study the $n$-th trace of the Millson theta form $\Theta_L^{\Millson}(\tau,Z)$ appearing in the theta lift.
	\item By splitting $\Theta_L^{\Millson}(\tau,Z)$ along $T_n$ into Siegel-type theta functions for sublattices of $L$, we explicitly compute $\tr_n(\Theta_L^{\Millson}(\tau,Z))$ in terms of unary theta functions (see Theorem~\ref{thm:EvaluationMillsonTX}). In total, we find that $\tr_n(\Omega_m^{\cusp})$ can be interpreted as a regularized Petersson inner product involving the harmonic Maass form $f_{m}$ and certain unary theta functions.
	\item An application of Stokes' Theorem in the spirit of the recent work of Bruinier, Ehlen, and Yang \cite{bruinier2020greens} yields a formula for $\tr_n(\Omega_m^{\cusp})$ which involves the Fourier coefficients of $f_{m}$, as well as the Fourier coefficients of holomorphic unary theta functions and their $\xi$-preimages. Recalling that $\tr_n(\Omega_m^{\cusp})$ is essentially equal to $\tr_m(\omega_{n}^{\mero})$, this gives the desired explicit formula for $\tr_m(\omega_{n}^{\mero})$. 
	\item Using results about the algebraicity of the Fourier coefficients of the harmonic Maass forms appearing in this formula, we derive the rationality of $\tr_m(\omega_n^{\mero})$.
\end{enumerate}

\subsection*{Outlook}
We expect that these results can be generalized to the more general setting of orthogonal modular forms for signature $(2, n)$. In particular, the $\xi$-operator can be defined in the same way, see also \cite{brfk} for case of the symplectic group $\Sp_{2g}$. Generalizations of Section \ref{sec:LocallyHarmonicMaassForms}, \ref{sec:ThetaLifts}, \ref{sec:SingRealAnalytic}, and \ref{subsec:CurrentEquation} then seem quite straightforward. However, the evaluation of regularized inner products for forms with singularities along algebraic cycles appears to be delicate. The authors plan to come back to this problem in the near future. Moreover, it would be interesting to drop the technical assumption, that the algebraic cycle $T_X$ and the real analytic cycle $C_Y$ do not intersect.

Integrals of the holomorphic functions $\omega_m(z_1, z_2)$ along the real analytic cycles have been considered by \cite{ZagierCycleIntegralsOmega, ZagierTracesHecke} in special cases. Zagier shows that these are related to the traces of Hecke operators and deduces a rationality result for these cycle integrals. The authors plan to investigate in the near future whether the methods used in this paper also lead to a generalization of these results.

\subsection*{Outline of the paper}
This work is organized as follows. In Section~\ref{sec:HarmonicMaassForms}, we review the basics of vector-valued harmonic Maass forms for the Weil representation. In Section~\ref{sec:OrthogonalModularForms} we set up the general framework, explaining the Grassmannian model of $\H^2$ and defining $\omega_{m}^{\cusp}$ and $\omega_n^{\mero}$ as orthogonal modular forms. In Section~\ref{sec:XiOperator} we introduce the $\xi$-operator on differential forms on $\Gamma \backslash \H^2$. In Section~\ref{sec:LocallyHarmonicMaassForms}, we examine the locally harmonic Hilbert-Maass forms $\Omega_m^{\cusp}$ and their relation to $\omega_m^{\cusp}$ via the $\xi$-operator.

In Section~\ref{sec:ThetaFunctions} we define the Millson theta form $\Theta_{L}^{\Millson}(\tau,Z)$ on $\H^2$ and evaluate it along the algebraic cycle $T_X$ using unary theta functions. In Section~\ref{sec:ThetaLifts}, we define the Millson theta lift $\Phi_L^{\Millson}(g,Z)$, showing that the locally harmonic Hilbert-Maass form $\Omega_m^{\cusp}$ is a lift of a harmonic Maass form. In Section~\ref{sec:InnerProducts} we evaluate regularized Petersson inner products between modular forms with singularities along real and algebraic cycles, proving $\omega_m^{\mero}$ is orthogonal to cusp forms and computing $\langle \omega_m^{\cusp},\omega_n^{\mero}\rangle_{\Pet}^{\reg}$ in terms of traces $\tr_m(\omega_n^{\mero})$ and $\tr_n(\Omega_m^{\cusp})$. We also prove the current equation for the Millson lift and its relation to the Oda lift.

Finally, in Section~\ref{sec:CycleIntegrals}, we derive an explicit formula for the traces of cycle integrals $\tr_m(\omega_n^{\mero})$ in terms of Fourier coefficients of harmonic Maass forms, presenting the general statement and proof of our main result, Theorem~\ref{thm:MainResultIntro}.

\section*{Acknowledgements}
We thank Jan Bruinier, Jens Funke, and Jan Vonk for helpful discussions.

\section{Harmonic weak Maass forms for the Weil representation}
\label{sec:HarmonicMaassForms}

In this section we recall some facts about the Weil representation and harmonic weak Maass forms. Throughout this section $L$ denotes an even lattice with a quadratic form $Q$ of signature $(b^+,b^-)$. The corresponding bilinear form is denoted by $(\cdot,\cdot)$, and $L'$ is the dual lattice.

\subsection{The Weil representation} We write $\Mp_2(\IR)$ for the metaplectic group, i.e., the double cover of $\SL_2(\IR)$ given by pairs $(M, \phi)$, where $M = \left(\begin{smallmatrix}a & b \\ c & d\end{smallmatrix}\right) \in \SL_2(\IR)$ and $\phi : \IH \to \IC$ is a holomorphic square-root of $c \tau + d$. 
The subgroup $\Mp_2(\IZ)$ is generated by the two elements $S = \left(\left(\begin{smallmatrix}0 & -1 \\ 1 & 0\end{smallmatrix}\right), \sqrt{\tau}\right)$ and $T = \left(\left(\begin{smallmatrix}1 & 1 \\ 0 & 1\end{smallmatrix}\right), 1\right).$

Let now $L$ be an even lattice of signature $(b^+, b^-)$ and consider the group ring $\IC[L' / L]$ with basis $(\frake_\beta)_{\beta \in L' / L}$ and inner product $\langle \frake_\beta, \frake_{\gamma} \rangle = \delta_{\beta, \gamma}$. This is anti-linear in the second variable. Moreover, we use the notation $\frake_\beta(x) = e(x) \frake_\beta = e^{2 \pi i x} \frake_\beta$.

There is a unitary representation $\rho_L : \Mp_2(\IZ) \to \GL(\IC[L' / L])$ called the Weil representation defined by
$$\rho_L(S) \frake_\beta = \frac{\sqrt{i}^{b^- - b^+}}{\sqrt{\lvert L' / L \rvert}} \sum_{\gamma \in L' / L} \frake_{\gamma}(-(\beta, \gamma)) \qquad \text{and} \qquad \rho_L(T) \frake_\beta = \frake_\beta(Q(\beta)).$$
We let $L^- = (L,-Q)$. Then the corresponding Weil representation $\rho_{L^-}$ is the dual Weil representation of $L$.

\subsection{Harmonic Maass forms} Let $k \in \frac{1}{2}\Z$. For functions $f : \IH \to \IC[L' / L]$ and $(M, \phi) \in \Mp_2(\IZ)$ we define the slash operator by
\[
(f \vert_{k, L} (M, \phi))(\tau) = \phi(\tau)^{-2k} \rho_L^{-1}(M, \phi) f(M \tau).
\]

Recall the Maass operators of weight $k$ on smooth functions $g : \IH \to \IC[L' / L]$. They act component-wise and are given by
$$L_k = -2 i v^2 \frac{\partial}{\partial \overline{\tau}} \qquad \text{and} \qquad R_k = 2i v^{-k} \frac{\partial}{\partial \tau} v^k = 2i \frac{\partial}{\partial \tau} + k v^{-1}, \qquad (\tau=u+iv).$$
They lower and raise the weight of a modular form by $2$, respectively.

We will also need the iterated raising operators defined by
$$R_k^j = R_{k + 2(j - 1)} \cdots R_{k + 2} R_{k}.$$
Following \cite[Section 3]{brfu04} we define a harmonic weak Maass form of weight $k$ with respect to the Weil representation $\rho_{L}$ as a real analytic function $f : \IH \to \IC[L' / L]$ such that $f \vert_{k, L} (M, \phi) = f$ for all $(M, \phi) \in \Mp_2(\IZ)$, $f$ is annihilated by the weight $k$ invariant Laplace operator $\Delta_{k} = -R_{k - 2} L_k$, and $f(\tau)$ has at most linear exponential growth as $v \to \infty$. We denote the space of these functions by $H_{k, L}$. The spaces of weakly holomorphic modular forms, holomorphic modular forms, and cusp forms, are denoted by $M_{k,L}^!, M_{k,L},$ and $S_{k,L}$, respectively.

Let 
\[
\xi_k f = v^{k-2} \overline{L_k f} = R_{-k} v^\kappa \overline{f}.
\]
Then $\xi_k$ defines an anti-linear map $H_{k, L} \to M_{2-k, L^-}^!$. The preimage of the space of cusp forms of weight $2-k$ under $\xi_k$ is denoted by $H_{k, L}^+$. A harmonic weak Maass form $g \in H_{k, L}^+$ has a decomposition $g = g^+ + g^-$, where $g^+$ is holomorphic with Fourier expansion
\[
g^+(\tau) = \sum_{\beta \in L' / L} \sum_{\substack{m \gg -\infty}} c_g^+(m, \beta) \frake_\beta(m \tau)
\]
and $g^-$ is non-holomorphic but exponentially decreasing as $v \to \infty$.

\subsection{Maass-Poincar\'{e} series}
\label{subsec:MaassPoincareSeries}

Let $k \in \frac{1}{2}\Z$ with $k < 0$, and denote by $M_{\mu, \nu}$ the usual $M$-Whittaker function (see \cite{abramowitzStegun}, equation 13.1.32).
We define, for $s \in \C$ and $v \in \R_{> 0}$,
\begin{equation}\label{eq:whittakerM}
	\calM_{k, s}(v) = v^{-\frac{k}{2}} M_{-\frac{k}{2}, s - \frac{1}{2}}(v).
\end{equation}
Following \cite{brhabil}, for $\beta \in L'/L$ and $m > 0$ we define the vector-valued Maass-Poincar\'e series
\[f_{m,\beta}(\tau, s) = \frac{1}{2\Gamma(2s)} \sum_{(M, \phi) \in \tilde{\Gamma}_\infty \backslash \Mp_2(\Z)} \left(\calM_{k, s}(4\pi m v)e(-mu) \frake_\beta\right) \vert_{k, L^-}(\gamma, \phi)(\tau), \]
where $\tilde{\Gamma}_\infty$ is the subgroup generated by $T$. The series converges absolutely for $\Re(s) > 1$, and at the special point $s = 1 - \frac{k}{2}$, the function
\[f_{m,\beta}(\tau) = f_{m,\beta}\left(\tau, 1 - \frac{k}{2}\right)\]
defines a harmonic Maass form in $H_{k, L^-}^+$ with principal part $e(-m\tau)(\frake_\beta + \frake_{-\beta}) + \frakc$ for some constant $\frakc \in \C[L'/L]$. In particular, every harmonic Maass form $f \in H_{k, L^-}^+$
can be written as a linear combination 
\[f(\tau) = \frac{1}{2} \sum_{\beta \in L'/L} \sum_{m >0} c_f^+( -m,\beta)f_{m,\beta}(\tau).\]

\subsection{Operators on vector-valued modular forms} We let $A_{k,L}$ be the space of $\C[L'/L]$-valued functions on $\H$ that transform like modular forms of weight $k$ for $\rho_L$. Let $K$ and $L$ be even lattices. For $n\in\N_0$ and functions $f \in A_{k,K}$ and $g \in A_{l,L}$ with $k,l \in \frac{1}{2} \Z$ we define the $n$-th \emph{Rankin--Cohen bracket}
\begin{align}\label{eq:RankinCohen}
[f,g]_n \coloneqq \sum_{\substack{r,s \geq 0\\ r+s=n}} (-1)^s \frac{\Gamma(k+n) \Gamma(l+n)}{\Gamma(s+1) \Gamma(k+n-s) \Gamma(r+1) \Gamma(l+n-r)}  f^{(r)} \otimes g^{(s)},
\end{align}
where the tensor product of two vector-valued functions $f = \sum_{\mu}f_\mu \e_\mu \in A_{k,K}$ and $g = \sum_{\nu}g_\nu \e_\nu \in A_{l,L}$ is defined by 
\[
f \otimes g \coloneqq \sum_{\mu , \nu} f_\mu g_\nu \e_{\mu + \nu} \in A_{k+l,K \oplus L},
\]
and the derivative is normalized by $f^{(r)} = \frac{1}{(2\pi i)^r}\frac{\partial^r}{\partial \tau^r}f$.
The proof of the following formula can be found in \cite[Proposition~3.6]{bruinier2020greens}.

\begin{Proposition}\label{rankin cohen brackets}
	Let $f \in H_{k,K}$ and $g \in H_{l,L}$ be harmonic Maass forms. For $n \in \N_0$ we have
	\begin{align*}
	&(-4\pi)^n L_{k+l+2n}\left([f,g]_n \right)= \frac{\Gamma(k+n) }{\Gamma(n+1) \Gamma(k)} L_k (f) \otimes R_l^n (g) +  \frac{(-1)^n \Gamma(l+n)}{  \Gamma(n+1) \Gamma(l)} R_k^n (f) \otimes L_l (g).
	\end{align*}
\end{Proposition}

If $K \subseteq L$ is a sublattice of finite index, then a vector-valued modular form $f \in A_{k, L}$ can be viewed as a modular form in $A_{k, K}$. This is achieved via the natural map
\begin{align}\label{eq:ArrowDown}
\begin{array}{ccl}
	A_{k, L} & \rightarrow & A_{k,K} \\
	f & \mapsto & f_K
\end{array} \qquad 
(f_K)_\mu =\begin{cases}
	f_{\overline{\mu}}, & \text{if } \mu \in L'/K, \\
	0, & \text{if } \mu \notin L'/K,
\end{cases}
\end{align}
where $\overline{\cdot}$ denotes the reduction $L'/K \rightarrow L'/L$.

On the other hand, for any $\overline{\mu} \in L'/L$, let $\mu$ be a fixed preimage of $\overline{\mu}$ in $L'/K$. Then we have a map
\begin{align}\label{eq:ArrowUp}
\begin{array}{ccl}
	A_{k, K} & \rightarrow & A_{k,L} \\
	g & \mapsto & g^L
\end{array} \qquad 
(g^L)_{\overline{\mu}} = \sum_{\alpha \in L/K} g_{\alpha + \mu}.
\end{align}
These maps are adjoint to each other with respect to the scalar products on $\IC[L'/L]$, respectively $\IC[K'/K]$. 
See \cite{scheithauer2015} and \cite{bruinier2009} for proofs and more details regarding these constructions.

\section{Orthogonal modular forms for $\Orth(2,2)$}
\label{sec:OrthogonalModularForms}

\subsection{The Grassmannian model of $\H^2$}\label{section grassmannian model}

We consider the four-dimensional $\R$-vector space
\begin{align*}
V(\IR) = \left\{ \begin{pmatrix}
a&\nu'\\\nu&b 
\end{pmatrix} \,:\, a, \nu, \nu', b\in\R \right\}
\end{align*}
together with the quadratic form 
\[
Q(X)=-\det(X)
\]
 and the corresponding bilinear form 
 \[
 (X,Y) = -\mathrm{tr}(X Y^*),
 \]
 where $\left(\begin{smallmatrix} a&b\\c&d\end{smallmatrix}\right)^*=\left(\begin{smallmatrix} d&-b\\-c&a\end{smallmatrix}\right)$. The real quadratic space $(V(\IR),Q)$ has signature $(2,2)$. The group $\SL_2(\R)^2$ acts isometrically on $V(\R)$ by
\[
(g_1,g_2).X = g_1 X g_2^t.
\]
Moreover, analysing the Clifford algebra yields $\mathrm{Spin}_{V(\IR)} \simeq \SL_2(\IR)^2$.

We consider the orthogonal basis
\[
e_{1} = \begin{pmatrix}-1 & 0 \\ 0 & 1 \end{pmatrix}, \quad e_{2} = \begin{pmatrix}0 & 1 \\ 1 & 0 \end{pmatrix}, \quad e_{3} = \begin{pmatrix} 1 & 0 \\ 0 & 1\end{pmatrix} ,\quad e_{4} =\begin{pmatrix}0 & -1 \\ 1 & 0 \end{pmatrix},
\]
of $V(\R)$. Their norms are $Q(e_{1}) = Q(e_{2}) = 1$ and $Q(e_{3}) = Q(e_{4}) = -1$. For $z = x+iy \in\H$ we let $g_{z} = \left( \begin{smallmatrix}\sqrt{y} & x/\sqrt{y} \\ 0 & 1/\sqrt{y}\end{smallmatrix}\right)$ such that $g_{z}i = z$. For $Z = (z_{1},z_{2}) \in \H^2$ we let $(g_{z_{1}},g_{z_{2}})$ act on the above basis and obtain  the orthonormal basis
\begin{align*}
X_{1}(Z) &= \frac{1}{\sqrt{y_{1}y_{2}}}\begin{pmatrix}x_{1}x_{2}-y_{1}y_{2} & x_{1}\\ x_{2} & 1  \end{pmatrix}, \qquad
X_{2}(Z) = \frac{1}{\sqrt{y_{1}y_{2}}}\begin{pmatrix}x_{1}y_{2}+x_{2}y_{1} & y_{1}\\ y_{2} & 0  \end{pmatrix}, \\
X_{3}(Z) &= \frac{1}{\sqrt{y_{1}y_{2}}}\begin{pmatrix}x_{1}x_{2}+y_{1}y_{2} & x_{1}\\ x_{2} & 1  \end{pmatrix}, \qquad
X_{4}(Z) = \frac{1}{\sqrt{y_{1}y_{2}}}\begin{pmatrix}x_{1}y_{2}-x_{2}y_{1} & -y_{1}\\ y_{2} & 0  \end{pmatrix},
\end{align*}
of $V(\R)$. We define the quantities
\begin{align*}
M(Z) &= \sqrt{y_{1}y_{2}}(X_{1}(Z)+iX_{2}(Z)) = \begin{pmatrix}z_{1}z_{2} & z_{1} \\ z_{2} & 1 \end{pmatrix} , \\
M(Z)^{\perp} &= \sqrt{y_{1}y_{2}}( X_{3}(Z)+iX_{4}(Z)) = \begin{pmatrix}\overline{z}_{1}z_{2} & \overline{z}_{1} \\ z_{2} & 1 \end{pmatrix}.
\end{align*}
For $\gamma \in \SL_{2}(\R)^2$ we have the transformation rules
\begin{align*}
M(\gamma Z) &= (c_1z_{1}+d_1)^{-1}(c_2z_{2}+d_2)^{-1}\gamma.M(Z), \\
M(\gamma Z)^{\perp} &= (c_1\overline{z}_{1}+d_1)^{-1}(c_2z_{2}+d_2)^{-1}\gamma.M(Z)^{\perp}.
\end{align*}

We let $\Gr(V)$ be the Grassmannian of positive definite planes in $V(\R)$. By abuse of notation, we also write $M(Z)$ for the positive plane in $\Gr(V)$ spanned by $X_{1}(Z),X_{2}(Z)$ and $M(Z)^\perp$ for the negative plane spanned by $X_3(Z),X_4(Z)$. The map $Z \mapsto M(Z)$ yields a bijection $\H^2 \cong \Gr(V)$ which is compatible with the corresponding actions of $\SL_2(\R)^2$.

For $Z = (z_1,z_2) \in \H^2$ and $X = \left(\begin{smallmatrix}a & \nu' \\ \nu & b \end{smallmatrix} \right) \in V(\R)$ we define the functions
\begin{align}\begin{split}\label{eq polynomials}
q_{Z}(X) &= (X,M(Z)) = -bz_{1}z_{2} + vz_{1} + v'z_{2} - a, \\
p_{Z}(X) &= \frac{1}{y_{1}}(X,M(Z)^{\perp}) = \frac{1}{y_{1}}(-b\overline{z}_{1}z_{2} + v\overline{z}_{1} + v'z_{2} - a),
\end{split}
\end{align}
which satisfy
\begin{align}\label{eq invariances}
\begin{split}
q_{\gamma Z}(X) &= (c_1z_{1}+d_1)^{-1}(c_2z_{2}+d_2)^{-1}q_{Z}(\gamma^{-1}.X), \\
p_{\gamma Z}(X) &= (c_1z_{1}+d_1)(c_2z_{2}+d_2)^{-1}p_{Z}(\gamma^{-1}.X),
\end{split}
\end{align}
for $\gamma \in \SL_2(\R)^2$.

Let $L \subset V(\R)$ be an even lattice such that $L \otimes \R = V(\R)$ and write $V = L \otimes \Q$. Let $\Orth^+(L)$ be the subgroup of $\SL_2(\IR)^2$ which maps $L$ to $L$. Let $L'$ be the dual lattice of $L$ and define the \emph{discriminant group} of $L$ by $L' / L$. It is a finite abelian group and the quadratic form $Q$ on $L$ induces a quadratic form $Q : L' / L \to \IQ / \IZ$ on $L' / L$. We have a map $\Orth^+(L) \to \Orth(L' / L)$ and we denote by $\Gamma(L)$ the subgroup of $\Orth^+(L)$ that acts trivially on $L' / L$. Throughout the paper we let $\Gamma$ be a commensurable subgroup of $\Gamma(L)$ for some even lattice $L \subset V(\R)$ with $L \otimes \IR = V(\IR)$.

We are mostly interested in the lattices
\begin{align}\label{number field lattice}
L = \left\{ \begin{pmatrix}a & \nu' \\ \nu & b \end{pmatrix} \, : \, a,b \in \Z, \, \nu \in \mathcal{O}_F \right\},
\end{align}
where $\mathcal{O}_F$ is the ring of integers in some real quadratic number field $F=\Q(\sqrt{D})$ of discriminant $D > 0$, and $\nu \to\nu'$ denotes the conjugation in $F$. In this case we have $\Gamma(L) = \SL_2(\mathcal{O}_F)$.

\subsection{Real analytic and algebraic cycles}

For $Y \in V(\R)$ with $Q(Y) > 0$ we consider the real analytic submanifold of $\H^2$ defined by
\[
C_Y = \{ Z \in \H^2 \,:\, (Y,M(Z)^\perp) = 0 \} = \{ Z \in \H^2  \,:\,  p_Z(Y) = 0 \}.
\]
It corresponds to the set of positive planes in $\Gr(V)$ containing $Y$, and hence to the Grassmannian of the real quadratic space $Y^\perp$, which has signature $(1,2)$. In particular, $C_Y$ has real dimension two.

Similarly, for $X \in V(\R)$ with $Q(X) < 0$ we consider the algebraic submanifold
\[
T_X = \{ Z \in \H^2  \,:\,  (X,M(Z))= 0 \} = \{ Z \in \H^2  \,:\,  q_Z(X) = 0 \},
\]
which can be identified with the Grassmannian of the signature $(2,1)$ quadratic space $X^\perp$.
 Hence, $T_X$ has complex dimension $1$. 

Now let $\Gamma$ be a subgroup of $\SL_2(\IR)^2$ which is commensurable with $\Gamma(L)$, and let $\Gamma_X$ be the stabilizer of $X$ in $\Gamma$. Then $\Gamma_X$ acts on $T_X$, and $\Gamma_X \backslash T_X$ defines an algebraic cycle in $\Gamma \backslash \H^2$. The following lemma gives a parametrization of $\Gamma_X \backslash T_X$. 

\begin{Lemma}\label{lemma parametrization}
	Let $X \in V(\R)$ with $Q(X) < 0$ and let $\gamma = (\gamma_1, \gamma_2) \in \SL_2(\R)^2$ with $\gamma^{-1}.X = \sqrt{|Q(X)|}\, e_4$. Then we have an isomorphism
	$$\IH \to T_X, \quad z \mapsto \gamma (z, z)$$
	and an isomorphism
	$$\SL_2(\IR) \to \SL_2(\IR)^2_X, \quad \sigma \mapsto \gamma (\sigma, \sigma) \gamma^{-1},$$
	where $\SL_2(\IR)^2_X$ is the stabilizer of $X$ in $\SL_2(\IR)^2$. Using this isomorphism, we can identify $\Gamma_X$ with a subgroup of $\SL_2(\IR)$ such that
	$$\Gamma_X \bs \IH \to \Gamma_X \bs T_X, \quad z \mapsto \gamma (z, z)$$
	is an isomorphism.
\end{Lemma}

Finally, we describe when two real analytic and algebraic cycles intersect. 
\begin{Lemma} \label{lemma intersection cycles}
Let $X, Y \in V(\R)$ with $Q(X) < 0$ and $Q(Y) > 0$. Then $T_X \cap C_Y \neq \emptyset$ if and only if $(X, Y) = 0$.
\end{Lemma}

\begin{proof}
The intersection $T_X \cap C_{Y}$ corresponds to the positive definite planes in $\Gr(V)$ that contain $Y$ and are orthogonal to $X$. This directly implies that $(X, Y) = 0$ if the intersection is non-empty. Conversely, if $(X, Y) = 0$, then there exists some positive definite plane in $X^\perp$ that contains $Y$ since $X^\perp$ has signature $(2, 1)$.
\end{proof}

\subsection{Orthogonal modular forms}
\label{subsec:OrhogonalModularForms}

Let $L \subset V(\IR)$ be an even lattice with $L \otimes \IR = V(\IR)$ and $\Gamma \subset \SL_2(\IR)^2$ be a subgroup which is commensurable with $\Gamma(L)$. We say that a function $f: \H^2 \to \C$ transforms like a modular form of weight $(k_1,k_2) \in \Z^2$ for $\Gamma$ if it satisfies
\[
f(\gamma Z) = (c_1 z_1 + d_1)^{k_1}(c_2 z_2 + d_2)^{k_2}f(Z)
\]
for all $\gamma = (\gamma_1,\gamma_2)=\left(\left(\begin{smallmatrix} a_1&b_1\\c_1&d_1\end{smallmatrix}\right), \left(\begin{smallmatrix} a_2&b_2\\c_2&d_2\end{smallmatrix}\right)\right) \in \Gamma$ and all $Z = (z_1,z_2) \in \H^2$. If $f$ is also holomorphic on $\H^2$, we call it a (holomorphic) modular form. If the lattice $L$ contains no $2$-dimensional isotropic subspace, then by the Koecher principle, a modular form is automatically holomorphic at the cusps; otherwise this is an additional assumption. If it vanishes at the cusps, it is called a cusp form. If the lattice $L$ is derived from a number field as in \eqref{number field lattice}, we call the corresponding modular forms Hilbert modular forms.

\subsubsection{Cusp forms} \label{sec:cuspforms}
We recall the definition of the Hilbert cusp forms which were first studied by Zagier \cite{zagierdoinaganuma}. Let $k \geq 4$ be even. For $Y \in V$ with $Q(Y) > 0$ we consider the function 
\[
\omega_Y^{\cusp}(Z) = \sum_{\gamma \in \Gamma_Y \backslash\Gamma}q_{Z}(\gamma^{-1}. Y)^{-k},
\]
where $\Gamma_Y$ denotes the stabilizer of $Y$ in $\Gamma$ and $q_Z(Y)$ is the function defined in \eqref{eq polynomials}. It follows from the rules \eqref{eq invariances} that $\omega_Y^{\cusp}(Z)$ transforms like a modular form of parallel weight $k$ for $\Gamma$. Moreover, it is easy to check that it is holomorphic on $\H^2$ and vanishes at the cusps, hence it defines a cusp form of weight $k$. For $\beta \in L'/L$ and $m > 0$ the group $\Gamma$ acts with finitely many orbits on the set $L_{m,\beta} = \{Y \in L+\beta \, :\, Q(Y) = m \}$, so we can define the finite averages
\begin{align}\label{eq average}
\omega_{m,\beta}^{\cusp}(Z) = \sum_{Y \in \Gamma \backslash L_{m,\beta}}\omega_Y^{\cusp}(Z) = \sum_{Y \in L_{m,\beta}}q_Z(Y)^{-k}.
\end{align}

\subsubsection{Meromorphic modular forms} \label{subsec:MeromorphicModularForms}
Mirroring the above construction, we choose $X \in V$ with $Q(X) < 0$ and consider the function
\[
\omega_X^{\mero}(Z) = \sum_{\gamma \in \Gamma_X\backslash\Gamma}q_{Z}(\gamma^{-1}. X)^{-k}.
\]
It defines a meromorphic modular form of weight $k$ for $\Gamma$ which vanishes at the cusps and has poles along the algebraic cycle $T_X$ and its $\Gamma$-translates in $\H^2$. Similarly as above, we define their averages for $\mu \in L'/L$ and $n < 0$ by 
\begin{equation} \label{equation average mero}
	\omega_{n,\mu}^{\mero}(Z) = \sum_{X \in \Gamma \backslash L_{n,\mu}}\omega_X^{\mero}(Z) = \sum_{X \in L_{n,\mu}}q_Z(Y)^{-k}.
\end{equation}

\begin{Remark}
If we choose the lattice $L$ coming from a number field $F$ as in \eqref{number field lattice} and the group $\Gamma = \SL_2(\mathcal{O}_F)$, and $n < 0$ is an integer, then the function $\sum_{\mu \in L'/L}\omega_{n/D,\mu}^{\mero}(Z)$ agrees with the meromorphic Hilbert modular form $\omega_{n}^{\mero}(Z)$ defined in the introduction.
\end{Remark}

\subsubsection{Restriction to algebraic cycles} It is well known that the restriction of a Hilbert modular form of weight $(k_1,k_2)$ to the diagonal yields an elliptic modular form of weight $k_1 + k_2$. More generally, we can define the restriction of an orthogonal modular form to an algebraic cycle $T_X$. To this end, we consider the slash operator of $\gamma = (\gamma_1,\gamma_2) \in \SL_2(\R)^2$ on functions $f: \H^2 \to \C$ given by
 \[
 (f|_{(k_1,k_2)}\gamma)(z_1,z_2) = j(\gamma_1,z_1)^{-k_1}j(\gamma_2,z_2)^{-k_2}f(\gamma_1 z_1,\gamma_2 z_2),
 \]
 where we put $j(\gamma_i,z_i) = c_i z_i + d_i$, $i=1,2$, as usual. Then we define the restriction to $T_X$ as follows.

\begin{Proposition}\label{prop:Restriction}
	 Let $X \in V(\R)$ with $Q(X) < 0$ and let $\gamma \in \SL_2(\R)^2$ with $\gamma^{-1}.X = \sqrt{|Q(X)|}\, e_4$. If $f: \H^2 \to \C$ transforms like a modular form of weight $(k_1,k_2)$ under $\Gamma$, then the \emph{restriction of $f$ to $T_X$,}
	 \[
	 f|_{T_X}(z) := \big(\left(f|_{(k_1,k_2)}\gamma\right)(z,z)\big)\big|_{k_1+k_2}\gamma_1^{-1} = j(\gamma_2 \gamma_1^{-1},z)^{-k_2}f(z,\gamma_2 \gamma_1^{-1}z)
	 \]
	 transforms like an elliptic modular form of weight $k_1 + k_2$ under $\Gamma_X$ and is independent of the choice of $\gamma$. Under the parametrization of Lemma~\ref{lemma parametrization} this coincides with the usual restriction.
\end{Proposition}

The proof is a simple application of Lemma~\ref{lemma parametrization}.

\section{The $\xi$-operator on Hilbert modular surfaces}
\label{sec:XiOperator}

In this section we give the definition of a differential operator on $\IH^2$ and its quotients $\Gamma \bs \IH^2$ that maps certain \textit{harmonic} Hilbert modular forms to Hilbert cusp forms. We follow the exposition of \cite[Chapter 3, §1]{FreitagHilbert} and refer to the book of Wells \cite{Wells} for further information on the differential geometry involved.

The $(1,1)$-form on $\H^2$ given by
\[
\omega = -\frac{i}{2} \partial \overline{\partial} \log(y_1 y_2) = d \mu(z_1) + d \mu(z_2)
\]
is $\SL_2(\R)^2$-invariant, where \[
d \mu(z_i) = \frac{d x_i \wedge d y_i}{y_i^2} = \frac{i}{2}\frac{d z_i \wedge d \overline{z}_i}{y_i^2}.
\]
The corresponding volume form is given by
\[
d\mu(Z) = \frac{1}{y_1^2 y_2^2} d x_1 \wedge d y_1 \wedge d x_2 \wedge d y_2 =  d \mu(z_1) \wedge d \mu(z_2).
\]
Let $\Gamma \subseteq \SL_2(\IR)^2$ be a commensurable subgroup of $\Gamma(L)$ for some even lattice $L \subset V(\R)$. Then the above differential forms define differential forms on $\Gamma \bs \IH^2$. We will write $\calA^{p, q}$ for the sheaf of $(p, q)$-forms on $\Gamma \bs \IH^2$ and for $U \subseteq \Gamma \bs \IH^2$ we denote its sections over $U$ by $\calA^{p, q}(U)$. For $\alpha, \beta \in \calA^{p, q}(U)$ the Hodge-$\overline{*}$-operator is defined by the equality
$$\alpha \wedge \overline{*} \beta = \langle \alpha, \beta \rangle d\mu(Z).$$
This defines an anti-linear map $\overline{*} : \calA^{p, q} \to \calA^{2 - p, 2 - q}$ with $\overline{*} d\mu(Z) = 1$. 

Consider the sheaf $\calL_k$ of modular forms of weight $(k,k)$ on $\Gamma \bs \IH^2$. For $U \subseteq \Gamma \bs \IH^2$ open its sections $\calL_{k}(U)$ are given by functions $F : \pi^{-1}(U) \to \IC$ that transform with weight $(k,k)$. Here, $\pi : \IH^2 \to \Gamma \bs \IH^2$ is the usual projection map. If $\Gamma$ is torsion-free, $\calL_{k}$ is a hermitian line bundle with hermitian metric given by the Petersson metric $F(z_1, z_2) \overline{G(z_1, z_2)} (y_1 y_2)^k$ on the fiber. The dual bundle $\calL_k^*$ can be identified with $\calL_{-k}$ using the hermitian metric and the map $F(z_1, z_2) \mapsto (y_1 y_2)^k \overline{F(z_1, z_2)}$ defines an anti-linear bundle isomorphism $\calL_k \to \calL_{-k}$. Using this map, we can define the Hodge-$\overline{*}$-operator on $\calA^{p, q}(\calL_k)$, the differential forms with values in $\calL_k$. For $U \subseteq \Gamma \bs \IH^2$ open and $\omega \in \calA^{p, q}(U), \, F \in \calL_k(U)$, we define
$$\overline{*}_k(\omega \otimes F) := (\overline{*} \omega) \otimes \left((y_1 y_2)^k \overline{F(Z)}\right) \in \calA^{2 - p, 2 - q}(U, \calL_{-k}).$$
This yields a Hodge-$\overline{*}$-operator $\overline{*}_k : \calA^{p, q}(\calL_k) \to \calA^{2 - p, 2 - q}(\calL_{-k})$. We let 
\[
\xi_{k, Z} = \overline{*}_{k} \overline{\partial},
\]
so that for modular forms $F$ of weight $k$, i.e., global sections of $\calL_k$, we have 
\begin{align}\label{eq:DeltaOnDifferentialForms}
\Delta_k F = \xi_{-k, Z} \xi_{k, Z} F.
\end{align}
In particular, for modular forms $F$ of weight $k$ we have $\xi_{k, Z} F \in \calA^{2, 1}(\Gamma \bs \IH^2, \calL_{-k})$ and for $\omega \in \calA^{2, 1}(\Gamma \bs \IH^2, \calL_{-k})$, we have \[
\xi_{-k, Z} \omega = \overline{*}_{-k} \overline{\partial} \omega \in \calL_{k}(\Gamma \bs \IH^2).
\] 
The following lemma gives the  $\xi_{-k, Z}$-operator on $(2, 1)$-forms explicitly in terms of the operators $\xi_{k,z_i} F = 2i y_i^{k} \overline{\frac{\partial }{\partial \overline{z}_i}F}$. We note that elements in $\calA^{2, 1}(\Gamma \bs \IH^2, \calL_{-k})$ can be written as
\begin{equation}\label{eq:defH}
H = H_1(Z) \, d z_1 \wedge d \mu(z_2) + H_2(Z) \, d \mu(z_1) \wedge d z_2
\end{equation}
with $H_1, H_2$ transforming like modular forms of weight $(2-k, -k), (-k, 2-k)$ in $Z$.

\begin{Lemma}\label{lem:XiOnH}
Let $H \in \calA^{2, 1}(\Gamma \bs \IH^2, \calL_{-k})$ be given as in \eqref{eq:defH}.
Then we have
$$\xi_{-k, Z} H = -y_2^{-k} \xi_{2-k,z_1} H_1(Z) - y_1^{-k} \xi_{2-k,z_2} H_2(Z).$$
\end{Lemma}

\begin{proof}
This follows by a direct calculation. 
\end{proof}

We will need the following reformulation of Stokes' Theorem.

\begin{Lemma}\label{lem:StokesForXi}
Let $U \subseteq V \subseteq \Gamma \bs \IH^2$ be open with $\overline{U} \subseteq V$. Let $F, G \in \calL_{k}(V)$ be cuspidal, i.e., they vanish at all cusps, and assume that $G = \xi_{-k, Z} H$ for some $H \in \calA^{2, 1}(V, \calL_{-k})$. Then
$$\int_{U} F(Z) \overline{G(Z)} (y_1 y_2)^k d\mu(Z) = -\int_{U} \xi_{k, Z} F(Z) \wedge \overline{*}_{-k} H(Z) + \int_{\partial U} F(Z) H(Z).$$
\end{Lemma}

\begin{proof}
We note that $\overline{G(Z)} (y_1 y_2)^k d\mu(Z) = \overline{*}_k G(Z)$ and obtain
\begin{align*}
\int_{U} F(Z) \overline{G(Z)} (y_1 y_2)^k d\mu(Z) = \int_U F(Z) \wedge \overline{*}_k \xi_{-k, Z} H(Z).
\end{align*}
Using $\overline{*}_k \overline{*}_{-k} = \id$ on $4$-forms yields $\overline{*}_k \xi_{-k, Z} H(Z)=\overline{\partial} H(Z)$.
Since $H \in \calA^{2, 1}(V, \calL_{-k})$ we have
$$d(F(Z) \wedge H(Z)) = \overline{\partial}(F(Z) \wedge H(Z)) = \overline{\partial} F(Z) \wedge H(Z) + F(Z) \wedge \overline{\partial} H(Z)$$
and thus
$$\int_U F(Z) \wedge \overline{\partial} H(Z) = - \int_U \overline{\partial} F(Z) \wedge H(Z) + \int_U d(F(Z) \wedge H(Z)).$$
Using the definition of $\overline{*}$ and that it is an isometry, we obtain
$$\overline{\partial} F(Z) \wedge H(Z) = \overline{*}_k \overline{\partial} F(Z) \wedge \overline{*}_{-k} H(Z).$$
Stokes' Theorem and the cuspidality of $F$ yield the result.
\end{proof}

\begin{Remark}
Let $\calL_{k, L}$ be the vector-bundle of modular forms of weight $k$ with respect to the Weil representation $\rho_L$.
The same procedure for $\SL_2(\IZ) \bs \IH$ and $\calL_{k, L}$ yields a differential operator $\xi_k : \calL_{k, L} \to \calA^{1, 0}(\calL_{-k, L^-}) \simeq \calL_{2-k, L^-}$ and a differential operator $\xi_{-k} : \calA^{1, 0}(\calL_{-k, L^-}) \simeq \calL_{2 - k, L^-} \to \calL_{k, L}$, where one uses $\calA^{1, 0} \simeq \calL_2$. This explains the notation $\xi_{2 - k}$ instead of $\xi_{-k}$ in \cite{brfu04}. 
In the Siegel case an analogous operator has been introduced in \cite{brfk}.
\end{Remark}

\section{Locally harmonic Hilbert-Maass forms}
\label{sec:LocallyHarmonicMaassForms}

In this section we define locally harmonic Hilbert-Maass forms as differential $(2,1)$-forms on $\Gamma \backslash \H^2$ and state their relation to the cusp forms $\omega_{Y}^{\cusp}(Z)$ defined in Section~\ref{sec:cuspforms}.

\begin{Definition}\label{def:locallyharmonic}
	For $Y \in V$ with $Q(Y) > 0$ and $Z \in \H^2$ not lying on any of the $\Gamma$-translates of the real analytic cycle $C_Y$ we define the \emph{locally harmonic Hilbert-Maass form}
	\[
	\Omega_{Y}^{\cusp}(Z) = \Omega_{Y, 1}^{\cusp}(Z) \, d z_1 \wedge d \mu(z_2) + \Omega_{Y, 2}^{\cusp}(Z) \, d \mu(z_1) \wedge d z_2,
	\]
	where
\begin{align*}
	&\Omega_{Y, 1}^{\cusp}(Z) = \sum_{\gamma \in \Gamma_Y \backslash \Gamma} y_1^{k - 2} y_2^k \frac{\overline{q_Z(\gamma^{-1}. Y)}^{1 - k}}{\overline{p_Z(\gamma^{-1}. Y)}}, \qquad
	\Omega_{Y, 2}^{\cusp}(Z) = \sum_{\gamma \in \Gamma_Y \backslash \Gamma} y_1^{k} y_2^{k - 2} \frac{\overline{q_Z(\gamma^{-1}. Y)}^{1 - k}}{\frac{y_1}{y_2}p_Z(\gamma^{-1}. Y)}.
\end{align*}
	For $\beta \in L'/L$ and $m >0$ and $Z \in \H^2$ not lying on any of the cycles $C_Y$ for $Y \in L_{m,\beta}$ we define the averages
	\[
	\Omega_{m,\beta}^{\cusp}(Y) = \sum_{Y \in \Gamma \backslash L_{m,\beta}}\Omega_{Y}^{\cusp}(Z). 
	\]
\end{Definition}
Using the invariances \eqref{eq invariances} we see that $\Omega_{Y, 1}^{\cusp}(Z)$ and $\Omega_{Y, 2}^{\cusp}(Z)$ transform like modular forms of weight $(2-k,-k)$ and $(-k,2-k)$ for $\Gamma$, respectively. Hence, $\Omega_{Y}^{\cusp}(Z)$  defines a $(2,1)$-form on $\Gamma \backslash \H^2$ with values in $\mathcal{L}_{-k}$ and singularities along the $\Gamma$-translates of the real analytic cycle $C_Y$. It is related to the cusp forms $\omega_Y^{\cusp}(Z)$ defined in Section~\ref{sec:cuspforms} by the $\xi$-operator on differential forms defined in Section~\ref{sec:XiOperator} as follows.

\begin{Proposition}\label{prop:XiOnLocallyHarmonic}
	For $Z \in \H^2$ outside of the singularities of $\Omega_{Y}^{\cusp}(Z)$ we have
	\[
	\xi_{-k,Z}\Omega_{Y}^{\cusp}(Z) = -2(k-1)\omega_{Y}^{\cusp}(Z).
	\]
\end{Proposition}

\begin{proof}
	By Lemma~\ref{lem:XiOnH} we have to show that
	\[
	\xi_{2-k,z_1}\Omega_{Y,1}^{\cusp}(Z) = (k-1)y_2^{k} \,\omega_{Y}^{\cusp}(Z), \qquad \xi_{2-k,z_2}\Omega_{Y}^{\cusp}(Z) = (k-1) y_1^{k}\, \omega_{Y}^{\cusp}(Z).
	\]
	This can be checked by a direct computation using that $y_1 \overline{p_Z(Y)}$ is holomorphic in $z_1$, $p_Z(Y)$ is holomorphic in $z_2$, and the rules
	\[
	\frac{\partial}{\partial \overline{z}_1}\left(\frac{1}{y_1}\overline{q_Z(Y)}\right) = -\frac{i}{2y_1}\overline{p_Z(Y)}, \qquad \frac{\partial}{\partial \overline{z}_2}\left(\frac{1}{y_2}\overline{q_Z(Y)}\right) = -\frac{iy_1}{2y_2^2}p_Z(Y).\qedhere
	\]
\end{proof}

Since $\omega_{Y}^{\cusp}(Z)$ is holomorphic on $\H^2$, we obtain:

\begin{Corollary}
	Outside of the singularities, the $(2,1)$-form $\Omega_Y^{\cusp}(Z)$ is harmonic with respect to the Laplace operator $\Delta_{-k} = \overline{*}_{k} \overline{\partial} \overline{*}_{-k} \overline{\partial} + \overline{\partial} \overline{*}_{k} \overline{\partial} \overline{*}_{-k} = \xi_k \xi_{-k} + \overline{\partial} \overline{*}_{k} \overline{\partial} \overline{*}_{-k}$. 
\end{Corollary}

\begin{Remark}
For $X \in V$ with $Q(X) > 0$ and $Z \in \H^2$ not lying on any of the $\Gamma$-translates of the algebraic cycle $T_X$ we define the \emph{polar harmonic Maass form} $\Omega_X^{\mero}(Z)$ by the same formulas as in Definition~\ref{def:locallyharmonic}. By the same computation as in the proof of Proposition~\ref{prop:XiOnLocallyHarmonic} the polar harmonic Maass form $\Omega_X^{\mero}(Z)$ is related to the meromorphic modular form $\omega_X^{\mero}(Z)$ by
\begin{align}\label{eq:OmegaXiPreimage}
	\xi_{-k, Z} \Omega_X^{\mero}(Z) = -2 (k-1) \omega_X^{\mero}(Z).
\end{align}
In particular, it is harmonic outside of the singularities.
\end{Remark}

\section{Evaluation of theta functions along algebraic cycles}
\label{sec:ThetaFunctions}

In this section we define the Millson theta form $\Theta_L^{\Millson}(\tau,Z)$ and compute its evaluation along an algebraic cycle $T_X$. We let $L \subset V(\IR)$ be an even lattice in the real quadratic space of signature $(2,2)$ defined in Section~\ref{section grassmannian model}, and we let $\Gamma \subset \SL_2(\R)^2$ be commensurable with $\Gamma(L)$. We let $V = L \otimes \Q$.

\subsection{Theta functions in signature $(2,2)$}\label{sec:O(2,2)theta}

 For brevity, we write $\lambda_Z$ and $\lambda_{Z^\perp}$ for the orthogonal projections of $\lambda \in V(\R)$ to the positive plane $M(Z)$ and the negative plane $M(Z)^\perp$, respectively. We consider the following Siegel-type theta functions.

\begin{Theorem}
\label{thm:thetafunctions} Let $k \in \Z$ with $k \geq 4$ be even. 
\begin{enumerate}
\item
The \emph{Doi-Naganuma theta function}
\[
\Theta_{L}^{\DN}(\tau,Z) = v\sum_{\lambda \in L'}\frac{q_{Z}(\lambda)^{k}}{(y_{1}y_{2})^{k}}e(Q(\lambda_{Z})\tau + Q(\lambda_{Z^{\perp}})\overline{\tau}) \e_{\lambda+L}
\]
has weight $k$ in $\tau$ for $\rho_{L}$, and $\overline{\Theta_L^{\DN}(\tau,Z)}$ has weight $(k,k)$ in $Z$ for $\Gamma$. 
\item
The \emph{Millson theta functions}
\begin{align*}
\Theta_{L,1}^{\Millson}(\tau,Z) &= v^{k}\sum_{\lambda \in L'} q_{Z}(\lambda)^{k-1}p_{Z}(\lambda) e(Q(\lambda_{Z})\tau + Q(\lambda_{Z^{\perp}})\overline{\tau}) \e_{\lambda+L}, \\
\Theta_{L,2}^{\Millson}(\tau,Z) &= v^{k} \frac{y_1}{y_2}\sum_{\lambda \in L'} q_{Z}(\lambda)^{k-1}\overline{p_{Z}(\lambda)} e(Q(\lambda_{Z})\tau + Q(\lambda_{Z^{\perp}})\overline{\tau}) \e_{\lambda+L}
\end{align*}
have weights $(2-k,-k)$ and $(-k,2-k)$ in $Z$ for $\Gamma$, respectively, and their complex conjugates have weight $2-k$ in $\tau$ for $\rho_{L^{-}}$.  
\end{enumerate}
\end{Theorem}

\begin{proof}
	For $k, l \in \N_0$ we consider the $\C[L'/L]$-valued Siegel theta function
\[
\Theta_{L}^{(l,k)}(\tau,Z) = \sum_{\lambda \in L'}q_{Z}(\lambda)^{l}p_{Z}(\lambda)^{k}e(Q(\lambda_{Z})\tau + Q(\lambda_{Z^{\perp}})\overline{\tau}) \e_{\lambda+L}.
\] 
Under the isomorphism $V \cong (\R^{4},x_{1}^{2}+x_{2}^{2}-x_{3}^{2}-x_{4}^{2})$ given by $x_i = (\lambda,X_i(Z))$ the function $q_{Z}(\lambda)^{l}p_{Z}(\lambda)^{k}$ corresponds to the polynomial $(x_{1}+ix_{2})^{l}(x_{3}+ix_{4})^{k}$, which is harmonic and homogeneous of degree $(l,k)$. Hence, \cite[Theorem~4.1]{borcherds} shows that the theta function transforms like a vector-valued modular form of weight $(1+l,1+k)$ for $\rho_L$ in $(\tau,\overline{\tau})$. Moreover, the invariances \eqref{eq invariances} imply that the theta function transforms like a modular form of weight $(-l+k,-l-k)$ in $Z$ under $\Gamma$. This implies the theorem.
\end{proof}

\begin{Definition}\label{def:millsontheta}
	We define the \emph{Millson theta form} as 
	\[
	\Theta_L^{\Millson}(\tau,Z) = \Theta_{L,1}^{\Millson}(\tau,Z)\, dz_1 \wedge d\mu(z_2) + \Theta_{L,2}^{\Millson}(\tau,Z)\, d\mu(z_1)\wedge dz_2.
	\]
	It defines a differential form in $\calA^{2, 1}(\Gamma \bs \IH^2, \calL_{-k})$.
\end{Definition}

 The Millson theta form and the Doi-Naganuma theta function are related by the following differential equation.
\begin{Proposition}\label{prop:difftheta}
We have
\begin{align*}
\xi_{-k, Z} \Theta_{L}^{\Millson}(\tau, Z) = 2\xi_{k, \tau} \Theta_L^{\DN}(\tau, Z),
\end{align*}
where $\xi_{-k,Z}$ is the $\xi$-operator on $(2,1)$-forms as defined in Section~\ref{sec:XiOperator}.
 
\end{Proposition}
\begin{proof}
By Lemma~\ref{lem:XiOnH} we need to show that
\begin{align*}
y_2^{-k} \xi_{2-k,z_{1}}\Theta_{L,1}^{\Millson}(\tau,Z)&= y_1^{-k} \xi_{2-k,z_{2}}\Theta_{L,2}^{\Millson}(\tau,Z) =  - \xi_{k,\tau}\Theta_L^{\DN}(\tau,Z).
\end{align*}
These equations can be checked using the identities
\begin{align*}
\notag Q(X_{Z}) &= \frac{1}{4y_{1}y_{2}}|q_{Z}(X)|^{2}, \qquad Q(X_{Z^{\perp}}) = -\frac{y_{1}}{4y_{2}}|p_{Z}(X)|^{2}, \\
\frac{\partial}{\partial \overline{z}_{1}}p_{Z}(X) &= -\frac{i}{2y_{1}^{2}}q_{Z}(X), \qquad \frac{\partial}{\partial \overline{z}_{1}}\frac{1}{y_{1}}\overline{q_{Z}(X)} = -\frac{i}{2y_{1}}\overline{p_{Z}(X)}.\qedhere
\end{align*}
\end{proof}

\begin{Remark}\label{remark theta arrow up} If $K \subset L$ is a sublattice of finite index, then the Doi-Naganuma theta functions for $K$ and $L$ are related via the map \eqref{eq:ArrowUp} by
\[
\Theta_L^{\DN}(\tau,Z) = \Theta_K^{\DN}(\tau,Z)^L.
\]
An analogous relation holds for the Millson theta form $\Theta_{L}^{\Millson}(\tau,Z)$.
\end{Remark}

\subsection{Theta functions in signature $(2, 1)$}
\label{subsec:ThetaFunctions21}

We now define two theta functions associated to lattices of signature $(2, 1)$. They will occur in a decomposition of the Millson theta function of signature $(2, 2)$ along algebraic cycles. Therefore, let $X \in V(\IR)$ and consider $W(\IR) = X^\perp$. Then $W(\IR)$ is a quadratic space of signature $(2, 1)$. For $z = x + i y \in \IH$ and $\lambda \in W(\IR)$ (or more generally, $\lambda \in V(\IR)$) set
$$q_z(\lambda) = q_Z(\lambda) \vert_{T_X}(z) \qquad \text{and} \qquad p_z(\lambda) = p_Z(\lambda) \vert_{T_X}(z).$$
Let $\gamma \in \SL_2(\IR)^2$ such that $\gamma^{-1}. X = \sqrt{|Q(X)|} \, e_4$. Then we write $\lambda_z = \lambda_{\gamma (z, z)}, \lambda_{z^\perp} = \lambda_{\gamma(z, z)^\perp}$ for $\lambda \in W(\IR)$. Moreover, the corresponding map $\SL_2(\IR) \to \SL_2(\IR)^2_X$ of Lemma~\ref{lemma parametrization} yields an action of $\SL_2(\IR)$ on $W(\IR)$. With this action we have
$$q_{\gamma z}(\lambda) = j(\gamma, z)^{-2} q_z(\gamma^{-1} \lambda) \qquad \text{and} \qquad p_{\gamma z}(\lambda) = p_z(\gamma^{-1} \lambda).$$

Let $P \subset W(\IR)$ be an even lattice such that $P \otimes \IR = W(\IR)$ and let $W = P \otimes \IQ$. Denote by $\Gamma(P)$ the subgroup of $\SL_2(\IR)$ that acts trivially on $P' / P$. 
\begin{Theorem}{\cite[Proposition 4.1]{alfesschwagenscheidtshintani}}
Let $k > 2$.
\begin{enumerate}
\item The \emph{Shintani theta function}
$$\Theta_P^{\Shintani}(\tau, z) = v^{k} \sum_{\lambda \in P'} q_z(\lambda)^{k-1} e(Q(\lambda_z) \tau + Q(\lambda_{z^\perp}) \overline{\tau}) \frake_{\lambda + P}$$
has weight $2-2k$ in $z \in \IH$ with respect to $\Gamma(P)$ and its complex conjugate has weight $1/2-k$ in $\tau\in\IH$ with respect to the Weil representation $\rho_{P^-}$.
 \item The  \emph{Millson theta function}
$$\Theta_P^{\Millson}(\tau, z) = v^{k} \sum_{\lambda \in P'} q_z(\lambda)^{k - 1} p_z(\lambda) e(Q(\lambda_z) \tau + Q(\lambda_{z^\perp}) \overline{\tau}) \frake_{\lambda + P}$$
 has weight $2-2k$ in $z \in \IH$ with respect to $\Gamma(P)$ and its complex conjugate has weight $3/2 - k$ in $\tau \in \IH$ with respect to the Weil representation $\rho_{P^-}$.
\end{enumerate}
\end{Theorem}

Let us assume that $P$ contains a primitive isotropic vector $\ell$. Choose $\ell' \in P'$ with $(\ell, \ell') = 1$. By replacing $\ell'$ with $\ell' - Q(\ell') \ell \in W$ if necessary we can assume that $\ell'$ is isotropic. Let
$$W_\ell = W\cap \ell^\perp \cap \ell'^\perp \qquad \text{and} \qquad K_\ell = W_\ell \cap P.$$
Then $W_\ell$ is positive definite one-dimensional rational quadratic space and $K_\ell \subset W_\ell$ is a positive definite even lattice. Let $w_\ell$ be a generator of $W_\ell(\R)$ with $(w_\ell,w_\ell) = 1$. We consider the unary theta function associated with $K_\ell$ and the polynomial $x^k$ on $W_\ell(\R)$ in the sense of \cite{borcherds}, that is,
\[
\Theta_{K_\ell, k+\frac{1}{2}}(\tau) = \sum_{\lambda \in K_\ell'} \exp\left(-\frac{\Delta}{8 \pi v}\right)(x^k)\big((\lambda, w_\ell)\big) e(Q(\lambda) \tau) \frake_{\lambda + K_\ell}.
\]
More explicitly, it is given by
\[
\Theta_{K_\ell,k+\frac{1}{2}}(\tau) = (2\sqrt{2\pi v})^{-k}\sum_{\lambda \in K_\ell'}H_{k}\left(\sqrt{2\pi v} (\lambda,w_\ell)\right) e(Q(\lambda)\tau) \e_{\lambda + K_\ell},
\]
with the Hermite polynomial $H_k(x) = (-1)^k e^{x^2} \frac{d^k}{dx^k}e^{-x^2}$. By \cite[Theorem~4.1]{borcherds}, the theta function $\Theta_{K_\ell,k+\frac{1}{2}}(\tau)$ transforms like a modular form of weight $k+1/2$ for $\rho_{K_\ell}$. For $k = 0$ it is a holomorphic modular form, and for $k = 1$ it is a cusp form. As explained in \cite[Section~2.2]{bif}, we can view $\Theta_{K_\ell,k+\frac{1}{2}}(\tau)$ as a modular form for $\rho_P$, and by a slight abuse of notation we will also denote this $\C[P'/P]$-valued function by $\Theta_{K_\ell,k+\frac{1}{2}}(\tau)$ in the following.

The unary theta functions of different weights are related by the raising operator as follows.

\begin{Lemma}\label{lem:RaisingUnaryThetas}
	We have $R_{k + \frac{1}{2}} \Theta_{K_\ell, k+\frac{1}{2}} = -2\pi\Theta_{K_\ell, k + \frac{5}{2}}$.
\end{Lemma}

The proof is a simple computation using some standard properties of Hermite polynomials, so we omit it for brevity.

By Witt's theorem, there exists $\sigma_\ell \in \SL_2(\IR)$ with
\[
\sigma_\ell^{-1} \ell = \begin{pmatrix}-1 & 0 \\ 0 & 0\end{pmatrix}, \qquad \sigma_\ell^{-1} \ell' = \begin{pmatrix}0 & 0 \\ 0 & 1\end{pmatrix}, \qquad \sigma_\ell^{-1} w_\ell = \frac{1}{\sqrt{2}}\begin{pmatrix}0 & 1 \\ 1 & 0 \end{pmatrix}.
\]
The next lemma gives us the asymptotics of the signature (2,1) Shintani and Millson theta functions as $z$ approaches the cusp of $\Gamma(P)$ corresponding to $\ell$.

\begin{Lemma}\label{lem:asymptotics}
	For $y \to \infty$ we have
	\begin{align*}
		\Theta_P^{\Shintani}(\tau, z) \vert_{2-2k} \sigma_\ell &= y^k v^{k-1/2} \big(\sqrt{2}i\big)^{k-1} \Theta_{K_\ell, k-\frac{1}{2}}(\tau)  + O(e^{-C y^2}), \\
		\Theta_P^{\Millson}(\tau, z) \vert_{2-2k} \sigma_\ell &= y^k v^{k - 3/2} \frac{ (k - 1) }{4\pi}\big(\sqrt{2}i\big)^{k - 2} \Theta_{K_\ell, k-\frac{3}{2}}(\tau) + O(e^{-C y^2}),
	\end{align*}
	for some constant $C > 0$.
\end{Lemma}

\begin{proof}
	In both cases one applies \cite[Theorem 5.2]{borcherds} and observes that there is only one term that yields a contribution. See also \cite[Section 5.4]{crawfordfunke}.
\end{proof}

\subsection{Evaluation of the Millson theta function along algebraic cycles}

Let $X \in V$ with $Q(X) < 0$. Then the lattices
$$P = L \cap (\IQ X)^\perp \subseteq W, \qquad N = L \cap \IQ X$$
have signature $(2, 1)$ and $(0, 1)$. Consider the theta functions
\begin{align}\label{eq:UnaryThetaFunctions}
\Theta_{N^-}(\tau) = \sum_{\lambda \in N'} e(-Q(\lambda) \tau) \frake_{\lambda + N}, \qquad \Theta^*_{N^-}(\tau) = \sum_{\lambda \in N'} \frac{(\lambda, X)}{\sqrt{\lvert Q(X) \rvert}} e(-Q(\lambda) \tau) \frake_{\lambda + N}.
\end{align}
Then $\Theta_{N^-}(\tau)$ is a holomorphic modular form of weight $1/2$ for $\rho_{N^-}$ and $\Theta_{N^-}^*(\tau)$ is a cusp form of weight $3/2$ for $\rho_{N^-}$.

\begin{Proposition}\label{prop:Splitting}
	For $z \in \IH$ we have the splitting
	$$\Theta_{L, 1}^{\Millson} \vert_{T_{X}}(\tau, z) = \left(\Theta_{P}^{\Millson}(\tau, z) \otimes \overline{\Theta_{N^-}(\tau)} + i \Theta_{P}^{\Shintani}(\tau, z) \otimes \overline{\Theta_{N^-}^{*}(\tau)}\right)^L.$$
	For $\Theta_{L, 2}^{\Millson} |_{T_{X}}(\tau, z)$ we have an analogous splitting, but with $i$ replaced by $-i$ in the second summand.
\end{Proposition}

\begin{proof}
	Using Remark \ref{remark theta arrow up} in this setting, we can reduce to the case $L = P \oplus N$. Then $L' = P' \oplus N'$ and every vector $\lambda \in L'$ can be uniquely written as a sum $\lambda = \lambda^P + \lambda^N$ for some $\lambda^P \in P', \lambda^N \in N'$. We have $Q(\lambda_z^N) = 0$ and $Q(\lambda_{z^\perp}^N) = Q(\lambda^N)$, so we obtain
	\begin{align*}
		&\Theta_{L, 1}^{\Millson} \vert_{T_X}(\tau, z) \\
		&= v^k \sum_{\substack{\lambda^P \in P' \\ \lambda^N \in N'}} q_{z}(\lambda^P + \lambda^N)^{k - 1} p_{z}(\lambda^P + \lambda^N) e(Q(\lambda_{z}^P) \tau + Q(\lambda_{z^\perp}^P) \overline{\tau}) e(Q(\lambda^N) \overline{\tau}) \frake_{\lambda^P} \otimes \frake_{\lambda^N}.
	\end{align*}
	We have
	$$q_{z}(\lambda^P + \lambda^N) = q_z(\lambda^P), \qquad p_{z}(\lambda^P + \lambda^N) = p_{z}(\lambda^P) + p_{z}(\lambda^N) = p_z(\lambda^P) + i \frac{(\lambda^N, X)}{\sqrt{\lvert Q(X) \rvert}}.$$
	This yields the result.
\end{proof}

 In order to define the ``evaluation at $T_X$'' of $\Theta_{L}^{\Millson}(\tau,Z)$, we introduce some notation, following \cite{funke}. Let us assume that $P$ contains an isotropic vector. We let $\Iso(P)$ be the set of isotropic lines in $P \otimes \Q$. The group $\Gamma_X$ acts on $\Iso(P)$ with finitely many orbits, and the set $\Gamma_X \backslash \Iso(P)$ is in bijection with the cusps of $\Gamma_X \backslash \H$. If we let $\Gamma_{X,\ell}$ denote the stabilizer of $\ell$ in $\Gamma_X$, then there is $\alpha_\ell \in \IR_{> 0}$ such that $\sigma_\ell^{-1} \Gamma_{X,\ell} \sigma_{\ell}$ is generated by the translations $\pm T^{\alpha_\ell}$ and we call $\alpha_\ell$ the width of the cusp $\ell$, see \cite[Definition 3.2]{funke}. For fixed $C \gg 0$ large enough, a fundamental domain for $\Gamma_X \backslash \H$ can be chosen as
	\[
	\mathcal{F}(\Gamma_X) = \calF_0 \cup \bigcup_{\ell \in \Gamma_X \backslash \Iso(P)}\sigma_\ell\mathcal{F}^{\alpha_\ell},
	\]
	where
	$$\calF^{\alpha_\ell} = \{ \tau = u + iv \in \IH : 0 \leq u \leq \alpha_\ell, v > C \}$$
	and $\calF_0$ is compact. We define the truncated fundamental domain for $T > C$ by replacing $\calF^{\alpha_\ell}$ with
	$$\calF^{\alpha_\ell}_T = \{ \tau = u + iv \in \IH : 0 \leq u \leq \alpha_\ell, T > v > C \}.$$

\begin{Definition}\label{def:EvaluationAtTX}
	Let $H(Z) = H_1(z) \, dz_1 \wedge d\mu(z_2) + H_2(z) \, d\mu(z_1) \wedge dz_2 \in \calA^{2, 1}(\Gamma \bs \IH^2, \calL_{-k})$. We define the \emph{evaluation of $H(Z)$ at $T_X$} by
	\[
	H(T_X) = \int_{\Gamma_X \backslash \H}^{\reg}R_{2-2k}^{k-1}\left(H_1|_{T_X}(z)-H_2|_{T_X}(z)\right)d\mu(z),
	\]
	with the restriction to $T_X$ from Proposition~\ref{prop:Restriction}, and the regularized integral defined by $\int_{\Gamma_X \backslash \H}^{\reg} = \lim_{T \to \infty}\int_{\mathcal{F}(\Gamma_X)_T}$.
\end{Definition}

Note that the integral is well-defined since the functions $H_i|_{T_X}(z)$ have weight $2-2k$ for $\Gamma_X$. We can now state the main result of this section.

\begin{Theorem}\label{thm:EvaluationMillsonTX}
	If $P$ is anisotropic, then $\Theta_L^{\Millson}(\tau,T_X) = 0$. Otherwise, we have
	\[
	\Theta_L^{\Millson}(\tau,T_X) =(-1)^{k/2+1} (k-2)!\sqrt{2}^{k+1}v^{k-1/2}\sum_{\ell \in \Gamma_X \backslash \Iso(P)} \alpha_\ell \left(\Theta_{K_\ell,k-\frac{1}{2}}(\tau) \otimes \overline{\Theta_{N^-}^*(\tau)} \right)^L.
	\]
\end{Theorem}

\begin{proof}
	Using Proposition~\ref{prop:Splitting} we see that
	\[
\Theta_L^{\Millson}(\tau,T_X) = \left(\int_{\Gamma_X \backslash \H}^{\reg}R_{2-2k,z}^{k-1}\left(2i\Theta_P^{\Shintani}(\tau,z) \otimes \overline{\Theta_{N^-}^*(\tau)}\right)d\mu(z)\right)^L.
\]
Hence, it suffices to compute the regularized integral over $\Gamma_X \backslash \H$ of $R_{2-2k,z}^{k-1}\Theta_P^{\Shintani}(\tau,z)$. 

	We split off the outermost raising operator and apply Stokes' Theorem in a similar form as in \cite[Lemma~4.2]{brhabil}. If $P$ is anisotropic, then $\Gamma_X \backslash \H$ is compact, so the integral vanishes. Otherwise, we compute
	\begin{align*}
	\int_{\Gamma_X \backslash \H}^{\reg}R_{2-2k,z}^{k-1}\Theta_P^{\Shintani}(\tau,z) d\mu(z) = -\sum_{\ell \in \Gamma_X\backslash \Iso(P)}\lim_{T \to \infty}\int_{0+iT}^{\alpha_\ell+iT} R_{2-2k,z}^{k-2}\left(\Theta_P^{\Shintani}(\tau,z)|_{2-2k}\sigma_\ell\right) \frac{dx}{y^2}.
	\end{align*}
	Now we plug in the asymptotics from Lemma~\ref{lem:asymptotics}. Since the square-exponentially decreasing parts hidden in the $O$-notation vanish as $T \to \infty$, we obtain
	\begin{align*}
	&\lim_{T \to \infty}\int_{0+iT}^{\alpha_\ell+iT} R_{2-2k,z}^{k-2}\Theta_P^{\Shintani}(\tau,z) \frac{dx}{y^2}  \\
	&\quad = \lim_{T \to \infty}\int_{0+iT}^{\alpha_\ell+iT} R_{2-2k,z}^{k-2}\left(y^k v^{k-1/2} \big(\sqrt{2}i\big)^{k-1} \Theta_{K_\ell, k-\frac{1}{2}}(\tau)\right)\frac{dx}{y^2}.
	\end{align*}
	We have $R_{k}y^{l} = (l + k)y^{l-1}$, so we inductively obtain $R_{2-2k}^{k-2}y^{k} = (-1)^{k-2}(k-2)!\, y^2$. Recall that $k$ is even, so $(-1)^{k-2} = 1$. Hence, the last expression simplifies to
	\begin{align*}
	v^{k-1/2}(k-2)!\big(\sqrt{2}i\big)^{k-1}\alpha_\ell\Theta_{K_\ell, k-\frac{1}{2}}(\tau).
	\end{align*}
	This finishes the proof.
\end{proof}

\begin{Remark}\label{rem:alphaell}
If we choose a generator $\kappa_\ell \in K_\ell$ with $w_\ell = \kappa_\ell / \lvert \kappa_\ell \rvert$, then $\kappa_\ell$ induces an \emph{Eichler transformation} $E(\kappa_\ell) \in \Gamma(P)_{\ell}$ and $\sigma_\ell^{-1} E(\kappa_\ell) \sigma_\ell$ acts via translation by $\sqrt{Q(\kappa_\ell)} = \sqrt{\lvert K_\ell' / K_\ell \rvert / 2}$. Since $\Gamma_{X}$ and $\Gamma(P)$ are commensurable, the cusp width $\alpha_\ell$ is a rational multiple of $\sqrt{\lvert K_\ell' / K_\ell \rvert / 2}$.
\end{Remark}

\section{The Millson theta lift}
\label{sec:ThetaLifts} 

Let $k \geq 4$ be an even integer. We recall the definition of the Doi-Naganuma theta lift of a cusp form $f \in S_{k,L}$ from Section~3.1 in Bruinier's part of \cite{123}. It is given by
\[
\Phi_L^{\DN}(f,Z) = \int_{\mathcal{F}}\left\langle f(\tau),\Theta_L^{\DN}(\tau,Z)\right\rangle v^{k} d\mu(\tau),
\]
where $\mathcal{F}$ is the standard fundamental domain for $\SL_2(\Z)\backslash\H$, $d\mu(\tau) = \frac{du \, dv}{v^2}$ is the invariant measure on $\H$, and $\Theta_L^{\DN}(\tau,Z)$ is the Doi-Naganuma theta function defined in Theorem~\ref{thm:thetafunctions}.

In this section we study a new theta lift which is related to the Doi-Naganuma lift via the $\xi$-operator, and which yields the locally harmonic Hilbert-Maass forms $\Omega_{m,\beta}^{\cusp}(Z)$ defined in Section~\ref{sec:LocallyHarmonicMaassForms}.

\begin{Definition} Let $g \in H_{2-k,L^-}^{+}$ be a harmonic Maass form of weight $2-k$ for the dual Weil representation $\rho_{L^-}$ which maps to a cusp form under $\xi_{2-k}$. We define the \emph{Millson theta lift} of $g$ by
\[
\Phi_L^{\Millson}(g,Z) = \int_{\mathcal{F}}^{\reg}\left\langle g(\tau),\overline{\Theta_{L}^{\Millson}(\tau,Z)}\right\rangle v^{2-k} d\mu(\tau),
\]
with the Millson theta form $\Theta_{L}^{\Millson}(\tau,Z)$ as in Definition~\ref{def:millsontheta}, and 
with the regularized integral defined by $\int_{\mathcal{F}}^{\reg} = \lim_{T \to \infty}\int_{\mathcal{F}_T}$, where $\mathcal{F}_T = \{\tau \in \H: v \leq T\}$ is a truncated fundamental domain.
\end{Definition}

	Note that, by the definition of the Millson theta form, the output of the Millson theta lift is a $(2,1)$-form on $\Gamma\backslash \H^2$. We first show that the Millson lift is well-defined outside of a certain set of singularities. To this end, we define a collection of real analytic cycles $C_Y$ corresponding to $g$ by
\[
C_g = \bigcup_{\substack{\beta \in L'/L,  m > 0 \\ c_g^+(-m,\beta) \neq 0}}\, \bigcup_{Y \in L_{m,\beta}}C_Y.
\]

\begin{Theorem}\label{thm:convergence}
The Millson theta lift $\Phi_L^{\Millson}(g,Z)$ converges point-wise for every $Z\in \mathbb{H}^2 \setminus C_g$.
\end{Theorem}

\begin{proof}
The proof is similar to \cite[Proposition~2.8]{brhabil}. We only treat the part of the Millson lift corresponding to $\Theta_{L,1}^{\Millson}(\tau, Z)$ since the proof for $\Theta_{L,2}^{\Millson}(\tau, Z)$ is analogous. First note that the integral 
\[
\int_{\mathcal{F}_1}\left\langle g(\tau),\overline{\Theta_{L,1}^{\Millson}(\tau,Z)}\right\rangle v^{2-k}d\mu(\tau)
\]
over the compact set $\mathcal{F}_1$ converges. Moreover, the non-holomorphic part $g^-$ of $g \in H_{2-k,L^{-}}^{+}$ is exponentially decreasing as $v\to\infty$. Since $\Theta_{L,1}^{\Millson}(\tau,Z)$ grows moderately, it is therefore enough to show the convergence of 
\begin{align*}
\int_{v=1}^\infty \int_{u=0}^1 \left\langle g^+(\tau),\overline{\Theta_{L,1}^{\Millson}(\tau,Z)}\right\rangle v^{2-k}d\mu(\tau).
\end{align*}
Plugging in the Fourier expansion of $g^+$ and writing out the theta function we get
\begin{align*}
\int_{v=1}^\infty &\int_{u=0}^1 \sum_{\beta\in L'/L} \sum_{n\gg-\infty } \sum_{\lambda\in L+\beta} c_g^+(n,\beta)  q_{Z}(\lambda)^{k-1}p_{Z}(\lambda)e(n\tau) e(Q(\lambda) u) e^{-2\pi v Q_Z(\lambda)} dudv,
\end{align*}
where $Q_Z(\lambda) = Q(\lambda_Z) - Q(\lambda_{Z^\perp})$. The integral over $u$
vanishes unless $n=-Q(\lambda)$, in which case it equals $1$. We obtain
\begin{equation}\label{eq:convergence1}
\int_{v=1}^\infty  \sum_{\beta\in L'/L}  \sum_{\lambda\in L+\beta} c_g^+(-Q(\lambda),\beta)  q_{Z}(\lambda)^{k-1}p_{Z}(\lambda) e^{2\pi v Q(\lambda)}e^{-2\pi v Q_Z(\lambda)} dv.
\end{equation}
We split the summation into three sums over $Q(\lambda)<0,\, Q(\lambda)=0, \, Q(\lambda)>0$, respectively. We only treat the case $Q(\lambda)>0$ here, since the other cases are similar, but simpler. The part of the integral in \eqref{eq:convergence1} for the terms satisfying $Q(\lambda)>0$ can be written as
\begin{equation}\label{eq:convergence2}
\sum_{\beta\in L'/L}  \sum_{m > 0} c_g^+(-m,\beta) \int_{v=1}^\infty 
\sum_{\substack{\lambda\in L+\beta\\Q(\lambda)=m}}   q_{Z}(\lambda)^{k-1}p_{Z}(\lambda) e^{-2\pi v m}e^{-2\pi v Q_Z(\lambda)} dv.
\end{equation}
Since $g^+$ has finite principal part, the first two sums are finite. Hence, it suffices to investigate the convergence of the integral for those finitely many $m > 0$ and $\beta \in L'/L$ for which $c_g^+(-m,\beta)\neq 0$. Note that the integral would clearly diverge if we had $Q_Z(\lambda) = m$ for some $\lambda \in L+\beta$ with $Q(\lambda)= m$. However, for $Q(\lambda) = m$ we have $Q_Z(\lambda) = m+\frac{y_1}{4y_2}|p_Z(\lambda)|^2$, so $Q_Z(\lambda) = m$ means $p_Z(\lambda) = 0$, or in other words, that $Z$ lies on the real analytic cycle $C_\lambda$. Since we excluded this case by taking $Z \in \H^2 \setminus C_g$, we can assume that $Q_Z(\lambda) \neq m$ for all $\lambda$ in the series.

As in the proof of \cite[Proposition~2.8]{brhabil} one can now reduce the convergence of the innermost series to the convergence of theta series for positive definite quadratic forms.
\end{proof}

We now show that the locally harmonic Hilbert-Maass form $\Omega_{m,\beta}^{\cusp}(Z)$ defined in Section~\ref{sec:LocallyHarmonicMaassForms} can be obtained as the Millson lift of the Maass-Poincar\'{e} series $f_{m,\beta}(\tau)$ of weight $2-k$ for $\rho_{L^-}$ defined in Section~\ref{subsec:MaassPoincareSeries}.

\begin{Theorem}\label{thm:unfoldpc}
For $Z$ outside of the singularities of the Millson theta lift, we have
\begin{align*}
\Phi_L^{\Millson}(f_{m,\beta},Z)
&= \frac{2}{\pi} (4m)^{k - 1} \Omega_{m, \beta}^{\cusp}(Z).
\end{align*}
\end{Theorem}
\begin{proof}
We compute the regularized integral of $f_{m,\beta}(\tau)$ against $\Theta_{L,1}^{\Millson}(\tau,Z)$. It is given by
\[
\frac{1}{2\Gamma(2s)}\int_{\mathcal{F}}^{\reg} \!\!\!\!\!\sum_{(M,\phi)\in \tilde{\Gamma}_\infty\setminus \Mp_2(\Z)} \left\langle  \left[\mathcal{M}_{2-k, s}(4\pi m v)\mathfrak{e}_\beta(-mu)\right]|_{2-k,L^-} (M,\phi),\overline{\Theta_{L,1}^{\Millson}(\tau,Z)}\right\rangle v^{2-k}d\mu(\tau).
\]
We follow the strategy of the proof of \cite[Theorem~2.14]{brhabil}. We apply the usual unfolding trick and insert the Fourier expansion of $\Theta_{L,1}^{\Millson}(\tau,Z)$ to obtain 
\begin{align*}
\frac{2}{\G(2s)} \int_{v=0}^\infty \int_{u=0}^1 \sum_{\lambda\in L+\beta} \mathcal{M}_{2-k, s}(4\pi m v) e(-mu) q_Z(\lambda)^{k-1} p_Z(\lambda) e(Q(\lambda)u)e^{-2\pi v Q_Z(\lambda)} dudv,
\end{align*}
where $Q_Z(\lambda) = Q(\lambda_Z) - Q(\lambda_{Z^\perp})$. The integral over $u$ vanishes unless $m=Q(\lambda)$, in which case it equals $1$. We plug in the definition of $\mathcal{M}_{2-k, s}(4\pi m v)$ from \eqref{eq:whittakerM} and arrive at
\begin{align*}
\frac{2}{\G(2s)} (4\pi m)^{k/2-1}& \sum_{\substack{\lambda\in L+\beta\\Q(\lambda)=m}}q_Z(\lambda)^{k-1} p_Z(\lambda) \int_{v=0}^\infty v^{k/2-1}  M_{k/2-1,s-1/2}(4\pi mv)e^{-2\pi v Q_Z(\lambda) } dv.
\end{align*}
The integral is a Laplace transform. By equation (11) on p.\@ 215 of \cite{laplacetransform} we have
\begin{align*}
\int_{v=0}^\infty v^{k/2-1}  M_{k/2 - 1,s-1/2}(4\pi mv)e^{-2\pi v Q_Z(\lambda)} dv &= (4\pi m)^s \G(s+k/2) \big(2\pi (m  + Q_Z(\lambda))\big)^{-s-k/2} \\
&\times{}_2F_{1}\left[s+k/2,s-k/2+1,2s;\frac{2m}{ m  + Q_Z(\lambda)}\right].
\end{align*}
Plugging in $s =  \frac{k}{2}$ we find that the regularized integral of $f_{m,\beta}(\tau)$ against $\Theta_{L,1}^{\Millson}(\tau,Z)$ equals
\begin{align*}
2 (4\pi m)^{k-1} &\sum_{\substack{\lambda\in L+\beta\\Q(\lambda)=m}} q_Z(\lambda)^{k-1} p_Z(\lambda) \big(2\pi (m  + Q_Z(\lambda))\big)^{-k} {}_2F_{1}\left[k,1,k;\frac{2m}{m + Q_Z(\lambda)}\right].
\end{align*}
We have ${}_2F_1[k,1,k;x]=\frac{1}{1-x}.$ Note that $Q_Z(\lambda)+m=2Q(\lambda_{Z})$ and $Q_Z(\lambda)-m = -2Q(\lambda_{Z^\perp})$ if $Q(\lambda) = m$. Therefore, we get
\[
-\frac{1}{2\pi}m^{k-1}\sum_{\substack{\lambda \in L+\beta\\Q(\lambda)=m}} q_Z(\lambda)^{k-1} p_Z(\lambda) \frac{Q(\lambda_{Z})^{1-k}}{Q(\lambda_{Z^\perp})}.
\]
If we now use $Q(\lambda_{Z}) = \frac{1}{4y_{1}y_{2}}|q_{Z}(\lambda)|^{2}$ and $Q(\lambda_{Z^{\perp}}) = -\frac{y_{1}}{4y_{2}}|p_{Z}(\lambda)|^{2}$, we obtain the statement of the theorem.
\end{proof}

We have the following relation between the Doi-Naganuma and the Millson theta lifts.
\begin{Theorem}\label{thm:ThetaLiftAndXiOperator}
Let $g\in H^+_{2-k,L^-}$. For $Z \in \H^2 \setminus C_g$ we have
\[
\xi_{-k,Z} \Phi_L^{\Millson}(g,Z) = -2 \Phi_L^{\DN}(\xi_{2-k,\tau} g, Z),
\]
or as a commutative diagram 
\begin{gather*}
\xymatrix@!=2.5cm {H^{+}_{2-k,L^-}\ar[r]^-{-2\xi_{2-k, \tau}}\ar[d]^-{\Phi_L^{\Millson}}  &\ar[d]^-{\Phi_L^{\DN}} S_{k,L} 
\\
\calA^{2, 1}(\Gamma \bs \IH^2, \calL_{-k}) \ar[r]^-{\xi_{-k,Z}} & \calL_k(\Gamma \bs \IH^2).}
\end{gather*}
\end{Theorem}
\begin{proof}
Using the relation $\xi_{-k,Z}\Theta_L^{\Millson}(\tau,Z) = 2\xi_{2-k,\tau}\Theta_L^{\DN}(\tau,Z)$ from Proposition~\ref{prop:difftheta} we get
\[
\xi_{-k,Z} \Phi_L^{\Millson}(g,Z)  = 2\int^{\reg}_\mathcal{F} \left\langle\xi_{2-k,\tau} \Theta_L^{\DN}(\tau,Z),g(\tau)\right\rangle v^{2-k} d\mu(\tau).
\]
By using Stokes' Theorem in the form as in \cite[Lemma 2.1]{brikavia} we obtain 
\[
\xi_{-k,Z} \Phi_L^{\Millson}(g,Z) = -2\Phi(\xi_{2-k,\tau}g,Z) -2\lim_{T\to\infty} \int_{\partial \mathcal{F}_T} \left\langle\overline{\Theta_L^{\DN}(\tau,Z)},g(\tau)\right\rangle d\bar\tau.
\]
The vanishing of the boundary integral can be shown using analogous arguments as in the proof of Theorem \ref{thm:convergence}.
\end{proof}

\section{Inner products of functions with singularities along real and algebraic cycles}
\label{sec:InnerProducts}

In this section we calculate regularized Petersson inner products of the form
\begin{align}\label{eq:preliminaryInnerProduct}
\left\langle G,\xi_{-k,Z}H \right\rangle_{\Pet}^{\reg} = \int_{\Gamma \bs \IH^2}^{\reg} G(Z) \overline{\xi_{-k,Z} H(Z)} (y_1 y_2)^k d\mu(Z),
\end{align}
where $G \in \calL_k(\Gamma \bs \IH^2)$ and $H \in \calA^{2,1}(\Gamma \bs \IH^2, \calL_{-k})$ are allowed to have certain singularities along algebraic cycles $T_X$ and real cycles $C_Y$, respectively. In particular, we will later apply our results to $G = \omega_{X}^{\mero}$ and $H = \Omega_{Y}^{\cusp}$. 

We first define suitable $\varepsilon$-neighbourhoods of the cycles $C_Y$ and $T_X$ that enable us to define the regularized inner product properly. 
Furthermore, to compute the regularized inner product effectively, certain integrals over the boundaries of the $\varepsilon$-neighbourhoods of $C_Y$ and $T_X$ must be calculated. These calculations will be carried out in the subsequent two subsections. 
In Section~\ref{subsec:RegularizedInnerProducts}, we combine our results on these integrals and give the explicit definition and a closed formula for the regularized inner product in~\eqref{eq:preliminaryInnerProduct}.

\subsection{Singularities along real analytic cycles}\label{sec:SingRealAnalytic}
Let $Y \in V$ with $Q(Y) > 0$. First assume $Y = e_2 = \left(\begin{smallmatrix}0 & 1 \\ 1 & 0 \end{smallmatrix} \right)$, so that the corresponding cycle is
$$C_Y = \{(z, - \overline{z}) \,:\,  z \in \IH \}$$
and the stabilizer of $Y$ is
$$G_Y = \{ (g, \tilde{g}) \,:\, g \in \SL_2(\IR) \}, \quad \widetilde{\begin{pmatrix}a & b \\ c & d\end{pmatrix}} = \begin{pmatrix}a & -b \\ -c & d\end{pmatrix}.$$
Write $k(\varphi) = \left(\begin{smallmatrix} \cos(2 \varphi) & \sin(2 \varphi) \\ -\sin(2 \varphi) & \cos(2\varphi) \end{smallmatrix}\right)$ and $a(e^{r}) = \left(\begin{smallmatrix} e^{r/2} & 0 \\ 0 & e^{-r/2} \end{smallmatrix}\right)$. Then we have a decomposition of $G = \SL_2(\R)^2$ given by
$$G = G_Y A_Y K,$$
where $K = \SO(2)^2$ and
$$A_Y = \{ (a(e^{-r}), a(e^r)) \in \SL_2(\IR)^2 : r \in \IR \}.$$
Let $\Gamma_Y = G_Y \cap \Gamma$. The map
$$\IH \times [0, \pi) \times [0, \infty) \to \IH \times \IH, \quad (z, \varphi, r) \mapsto (g_z k(\varphi) a(e^{-r}) i, \tilde{g}_z k(-\varphi) a(e^r) i)$$
is a diffeomorphism away from the cycle $C_Y$ and the set
$$B_\varepsilon(C_Y) = \{(g_z k(\varphi) a(e^{-r}) i, \tilde{g}_z k(-\varphi)) a(e^r) i \,:\, 0 \leq r < \varepsilon, \varphi \in [0, \pi]\}$$
is an open neighbourhood of the cycle $C_Y$. For arbitrary $Y$ choose $\gamma \in \SL_2(\IR)^2$ with $\gamma^{-1}.Y = \sqrt{Q(Y)} \, e_2$. Then $C_Y = \gamma C_{e_2}$ and we define $B_{\varepsilon}(C_Y) = \gamma B_{\varepsilon}(C_{e_2})$.

\begin{Lemma}\label{lem:RealCycleIntegrals}
Let $Y \in V$ with $Q(Y) > 0$ and let $G : \IH^2 \to \IC$ be cuspidal of weight $k$ for $\Gamma$ and smooth in $B_\varepsilon(C_Y)$. Then
\begin{align*}
&\lim_{\varepsilon \to 0} \int_{\Gamma_Y \bs \partial B_{\varepsilon}(C_Y)} G(Z) \left(y_1^{k-2}y_2^{k}\frac{\overline{q_Z(Y)}^{1 - k}}{\overline{p_Z(Y)}} \, d z_1 \wedge d \mu(z_2) + y_1^{k}y_2^{k-2}\frac{\overline{q_Z(Y)}^{1-k}}{\frac{y_1}{y_2}p_{Z}(Y)} \, d \mu(z_1) \wedge d z_2 \right) \\
&= -4\pi i \big(4Q(Y)\big)^{1-k} \int_{\Gamma_Y \bs C_Y} G(Z) q_Z(Y)^{k - 2} d z_1 \wedge d z_2.
\end{align*}
\end{Lemma}

\begin{proof}
First assume $Y = e_2$. Then a short calculation shows
\begin{align*}
	g_z k(\varphi) a(e^{-r})i &= x + \frac{y \sin(2 \varphi) \sinh(r)}{\cos(2 \varphi) \sinh(r) + \cosh(r)} + i \frac{y}{\cos(2\varphi) \sinh(r) + \cosh(r)}, \\
	\tilde{g}_z k(-\varphi) a(e^{r})i &= - x + \frac{y \sin(2 \varphi) \sinh(r)}{-\cos(2 \varphi) \sinh(r) + \cosh(r)} + i \frac{y}{-\cos(2\varphi) \sinh(r) + \cosh(r)}.
\end{align*}
In these coordinates, we have
\begin{align*}
	p_{z, \varphi, r}(e_2) &= p_{(g_z k(\varphi) a(e^{-r}) i, \tilde{g}_z k(-\varphi) a(e^{r}) i)}(e_2) \\
	&= 2i j(g_z k(\varphi) a(e^{-r}), i) j(\tilde{g}_z k(-\varphi) a(e^{r}) i, i)^{-1} \sinh(r).
\end{align*}
Moreover, on $\partial B_\varepsilon(C_{e_2})$ we obtain
\begin{align*}
	d z_1
	&= d x + \frac{\sin(2\varphi) \sinh(\varepsilon) + i}{\cos(2 \varphi) \sinh(\varepsilon) + \cosh(\varepsilon)} d y \\
	&+ 2 y \sinh(\varepsilon) \frac{\sinh(\varepsilon) + \cos(2 \varphi) \cosh(\varepsilon) + i \sin(2 \varphi)}{(\cos(2 \varphi) \sinh(\varepsilon) + \cosh(\varepsilon)))^2} d \varphi, \\
	d z_2
	&= - d x + \frac{-\sin(2\varphi) \sinh(\varepsilon) + i}{-\cos(2 \varphi) \sinh(\varepsilon) + \cosh(\varepsilon)} d y \\
	&+ 2 y \sinh(\varepsilon) \frac{-\sinh(\varepsilon) + \cos(2 \varphi) \cosh(\varepsilon) - i \sin(2 \varphi)}{(-\cos(2 \varphi) \sinh(\varepsilon) + \cosh(\varepsilon)))^2} d \varphi,
\end{align*}
so that
\begin{align*}
	\lim_{\varepsilon \to 0} \frac{d z_1 \wedge d z_2 \wedge d \overline{z}_2}{\overline{p_{z, \varphi, \varepsilon}(e_2)}}
	&= -2i y d z \wedge d \overline{z} \wedge d \varphi, \\
	\lim_{\varepsilon \to 0} \frac{d z_1 \wedge d \overline{z}_1 \wedge d z_2}{p_{z, \varphi, \varepsilon}(e_2)}
	&= -2i y d z \wedge d \overline{z} \wedge d \varphi.
\end{align*}
Recalling that $d\mu(z_i) = \frac{dx_i \wedge dy_i}{y_i^2} = \frac{i}{2}\frac{dz_i \wedge d\overline{z}_i}{y_i^2}$, we find that the expression on the left-hand side in the statement of the lemma equals
\begin{align*}
	2\int_{\Gamma_{e_2} \bs\IH} G(z, -\overline{z}) y^{2k-3}\overline{q_{(z, -\overline{z})}(e_2)}^{1-k} \left(\int_{0}^\pi d \varphi \right) dz \wedge d\overline{z}.
\end{align*}
Using $q_{(z,-\overline{z})}(e_2)=2iy$ and evaluating the inner integral to $\pi$ gives
\[
-4^{2-k}\pi i \int_{\Gamma_{e_2} \bs C_{e_2}} G(Z) q_{Z}(e_2)^{k-2} dz_1 \wedge d z_2.
\]
We can reduce the case of arbitrary $Y$ to the one above by choosing  $\gamma \in \SL_2(\IR)^2$ with $\gamma^{-1}.Y =  \sqrt{Q(Y)}\, e_2$ and using the transformation behaviour of the integrand.
\end{proof}

\subsection{Singularities along algebraic cycles}
Let $X \in V$ with $Q(X) < 0$. First suppose that $X = e_4 = \left( \begin{smallmatrix}0 & -1 \\ 1 & 0 \end{smallmatrix}\right)$, so that
$$T_{X} = \{(z, z) \in \IH^2 : z \in \IH \}$$
is the diagonal. We can identify $\IH^2$ with $\IH \times \ID$ using the real analytic map
$$(z_1, z_2) \mapsto \left(z_1, \frac{z_2 - z_1}{z_2 - \overline{z}_1}\right),$$
which maps the cycle $T_X$ to $\IH \times \{0\}$ with inverse
$$(z_1, w_2) \mapsto \left(z_1, \frac{z_1 - \overline{z}_1 w_2}{1 - w_2}\right).$$
For $\varepsilon > 0$ we write $B_\varepsilon(T_X)$ for the preimage of $\IH \times B_\varepsilon(0)$ under this map. This yields a neighbourhood $B_\varepsilon(T_X)$ of $T_X$. For general $X$ choose an element $\gamma \in \SL_2(\IR)^2$ such that $\gamma^{-1}.X = \sqrt{|Q(X)|}\, e_4$. Then $B_{\varepsilon}(T_X) = \gamma^{-1} B_{\varepsilon}(T_{e_4})$.

Our goal in this subsection is to evaluate integrals of the form 
\begin{equation}\label{eq:intalgebraiccycle}
 \lim_{\varepsilon \to 0}\int_{\Gamma_X\setminus\partial B_\varepsilon(T_X)} \frac{H(Z)}{q_Z(X)^k}
\end{equation}
for $H\in \mathcal{A}^{2,1}(\Gamma\setminus \mathbb{H}^2,\mathcal{L}_{-k})$. These calculations are far more involved than for the real cycles since the corresponding singularities have higher order. We inductively identify a primitive of $ \frac{H(Z)}{q_Z(X)^k}$ by defining an iterated raising operator on Hilbert modular forms, see Lemma \ref{lem:DifferentialOperatorCalculations}. On the cycle $T_{e_4}$ we can explicitly evaluate this expression (compare Lemma \ref{lem:DiagonalRaisingOperator}) which enables us to calculate the integral \eqref{eq:intalgebraiccycle} in Theorem \ref{thm:AlgebraicCycleIntegrals}. Our main result here can be viewed as a generalization of the Residue Theorem for $(2,1)$-forms valued in the sheaf of Hilbert modular forms of parallel weight $(-k,-k)$ vastly generalizing \cite[Lemma 4.1]{anbs}.

\begin{Lemma}\label{lem:intrealana}
	Let $f : B_\varepsilon(0) \to \IC$ be a real analytic function and $k\in\mathbb{N}$. Then
	$$\lim_{\varepsilon \to 0} \int_{\partial B_\varepsilon(0)} \frac{f(w)}{w^k} d \overline{w} = 0\quad\quad
	\text{and}\quad\quad
\lim_{\varepsilon \to 0} \int_{\partial B_\varepsilon(0)} \frac{f(w)}{w} d w = 2 \pi i f(0).$$
\end{Lemma}

\begin{proof}
	Since $f$ is real analytic, it has an expansion
	$$f(w) = \sum_{m \in \IZ} a_{m}(\lvert w \rvert) w^m, \quad a_m(\lvert w \rvert) = \sum_{n \geq 0} b_{m, n} \lvert w \rvert^n.$$
	In particular, $a_m(\lvert w \rvert)$ is bounded for $\lvert w \rvert \to 0$. Hence, using the residue theorem we obtain
	\begin{align*}
		\int_{\partial B_\varepsilon(0)} \frac{f(w)}{w^k} d \overline{w}
		&= \sum_{m \in \IZ} a_{m}(\varepsilon) \int_{\partial B_\varepsilon(0)} w^{m-k} d \overline{w} \\
		&= \sum_{m \in \IZ} a_{m}(\varepsilon) \varepsilon^{2m - 2k} \int_{\partial B_\varepsilon(0)} \overline{w}^{k-m} d \overline{w} \\
		&= (-2 \pi i) a_{k+1}(\varepsilon) \varepsilon^{2} \overset{\varepsilon \to  0}{\longrightarrow} 0.
	\end{align*}
	The second assertion is well-known and actually holds for continuous $f$.
\end{proof}

\begin{Corollary}\label{cor:localCalculations}
	Let $H(Z) d \mu(z_1) \wedge d \overline{z}_2 \in \calA^{2, 1}(\Gamma_X \bs \IH^2, \calL_{-k})$ be real analytic. For $k \in \IN$ we have
	\begin{align}\label{eq:VanishingIntegral}
		\lim_{\varepsilon \to 0} \int_{\Gamma_X \bs \partial B_{\varepsilon}(T_X)} \frac{H(Z)}{q_Z(X)^k} d\mu(z_1) \wedge d \overline{z}_2 = 0.
	\end{align}
	Moreover, 
	$$\lim_{\varepsilon \to 0} \int_{\Gamma_X \bs \partial B_{\varepsilon}(T_X)} \frac{H(Z)}{q_Z(X)} d \mu(z_1) \wedge d z_2 = -\frac{2 \pi i}{\sqrt{|Q(X)|}} \int_{\Gamma_X \bs T_X} H(Z) d \mu(z_1).$$
	Analogously, if $H(Z) d\overline{z}_1 \wedge d \mu(z_2) \in \calA^{2, 1}(\Gamma_X \bs \IH^2, \calL_{-k})$ and $k \in \IN$, then
	$$\lim_{\varepsilon \to 0} \int_{\Gamma_X \bs \partial B_{\varepsilon}(T_X)} \frac{H(Z)}{q_Z(X)^k} d\overline{z}_1 \wedge d \mu(z_2) = 0$$
	and
	$$\lim_{\varepsilon \to 0} \int_{\Gamma_X \bs \partial B_{\varepsilon}(T_X)} \frac{H(Z)}{q_Z(X)} dz_1 \wedge d \mu(z_2) = \frac{2 \pi i}{\sqrt{|Q(X)|}} \int_{\Gamma_X \bs T_X} H(Z) d \mu(z_2).$$
\end{Corollary}

\begin{proof}
	First assume $X = e_4$. Then $q_Z(e_4) = (z_1 - z_2)$. In the coordinates of $\IH \times \ID$, we have
	$$q_{z_1, w_2}(e_4) = -\frac{2i y_1}{1 - w_2} w_2 \quad \text{and} \quad d z_2 = \frac{1}{1 - w_2} d z_1 - \frac{w_2}{1 - w_2} d \overline{z}_1 + \frac{2 i y_1}{(1 - w_2)^2} d w_2.$$
	Hence, the integral in \eqref{eq:VanishingIntegral} is locally given by
	\begin{align*}
		\int_{z_1 \in U} \int_{\partial B_{\varepsilon}(0)} \frac{\tilde{H}(z_1, w_2)}{w_2^k} d \mu(z_1) \wedge d \overline{w}_2 
		= \int_{z_1 \in U} \left(\int_{\partial B_{\varepsilon}(0)} \frac{\tilde{H}(z_1, w_2)}{w_2^k} d \overline{w}_2\right) d \mu(z_1),
	\end{align*}
	where $\tilde{H}$ is a real analytic function. By Lemma \ref{lem:intrealana}, the inner integral vanishes for $\varepsilon \to 0$. The case of general $X$ can be reduced to this case.
	
	In the second case, the integral is locally given by
	\begin{align*}
		& -\lim_{\varepsilon \to 0} \int_{z_1 \in U} \left(\int_{\partial B_\varepsilon(0)} \frac{H(z_1, w_2)}{(1 - w_2) w_2} d w_2\right) d \mu(z_1) = -2 \pi i \int_{z_1 \in U} H(z_1, 0) d \mu(z_1).
	\end{align*}
	This shows the result for $X = e_4$ and the general case can again be reduced to this one.
\end{proof}

For differential forms with values in the line bundle $\calL_{-k}$ of modular forms of weight $-k$ we define the differential operator
$$\partial_{-k} = 2i \partial - \frac{k}{y_1} d z_1 - \frac{k}{y_2} d z_2 = R_{-k, z_1} d z_1 + R_{-k, z_2} d z_2 : \calA^{p, q}(\calL_{-k}) \to \calA^{p+1, q}(\calL_{-k}).$$
If
$H(Z) = H_1(Z) d z_1 \wedge d \mu(z_2) \in \calA^{2, 1}(\Gamma_X \bs \IH^2, \calL_{-k})$
we rewrite it as
$$H(Z) = \frac{1}{2} (\partial_{-1} q_Z(X)) \wedge \tilde{H}_1(Z) - \frac{1}{2} (\partial_{-1} q_Z(X)) \wedge \tilde{H}_2(Z)$$
with
\begin{align}\label{eq:tildeh1hand2}
	\tilde{H}_1(Z) &= \frac{H_1(Z)}{R_{-1, z_1} q_Z(X)} d \mu(z_2) \in \calA^{1, 1}(\Gamma_X \bs \IH^2, \calL_{1-k}), \\
	 \tilde{H}_2(Z) &= \frac{y_2^{-2} H_1(Z)}{R_{-1, z_2} q_Z(X)} d z_1 \wedge d \overline{z}_2 \in \calA^{1, 1}(\Gamma_X \bs \IH^2, \calL_{1-k}).
\end{align}
Define $\tilde{\partial}_{1-k} H(Z) = \partial_{1-k} (\tilde{H}_1(Z) + \tilde{H}_2(Z)) \in \calA^{(2,1)}(\Gamma_X \bs \IH^2, \calL_{1-k})$. The operator $\tilde{\partial}_{1-k}$ depends on the choice of $X$. We define
$$\tilde{\partial}_{1 - k}^j = \tilde{\partial}_{1 - k + j} \cdots \tilde{\partial}_{2 - k} \tilde{\partial}_{1 - k}.$$

\begin{Lemma}\label{lem:DifferentialOperatorCalculations}
	We have
	\begin{align*}
		\tilde{\partial}_{1-k}^j H(Z)
		&= \sum_{i = 0}^j \binom{j}{i} \frac{R_{2 - k, z_1}^{j - i} R_{-k, z_2}^i H_1(Z)}{(R_{-1, z_1} q_Z(X))^{j - i} (-R_{-1, z_2} q_Z(X))^i} d z_1 \wedge d \mu(z_2) \\
		&\quad + G(Z) (R_{-1, z_1} R_{-1, z_2} q_Z(X)) d z_1 \wedge d \mu(z_2) \in \calA^{2, 1}(\Gamma_X \bs \IH^2, \calL_{j-k}),
	\end{align*}
	where $G(Z)$ is a smooth function that has no singularities along $T_X$.
\end{Lemma}

\begin{proof}
	We proceed by induction on $j$. Using that  $R_{-1, z_i}^2 q_Z(X) = 0$ and the product rule for $R_{k, z_i}$, e.g., for functions $f,g:\mathbb{H}\to\mathbb{C}$ we have
	\begin{equation}\label{eq:productrarising}
	R_{k_1 + k_2, z_i} (f g) = (R_{k_1, z_i} f) g + f (R_{k_2, z_i} g),
	\end{equation}
	yields
	\begin{align*}
		\tilde{\partial}_{1 + k} H(Z) 
		&= \left(R_{1 + k, z_1} \frac{H_1(Z)}{R_{-1, z_1} q_Z(X)} - y_2^2 R_{1 + k, z_2} \frac{y_2^{-2} H_1(Z)}{R_{-1, z_2} q_Z(X)}\right) d z_1 \wedge d \mu(z_2) \\
		&= \left(\frac{R_{2 + k, z_1} H_1(Z)}{R_{-1, z_1} q_z(X)} - \frac{R_{k, z_2} H_1(Z)}{R_{-1, z_2} q_z(X)}\right) d z_1 \wedge d \mu(z_2) \\
		&= \sum_{i = 0}^1 \binom{1}{i} \frac{R_{2 + k, z_1}^{1 - i} R_{k, z_2}^{i} H_1(Z)}{(R_{-1, z_1} q_Z(X))^{1 - i} (-R_{-1, z_2} q_Z(X))^i} d z_1 \wedge d \mu(z_2).
	\end{align*}
	If the assertion holds for $j$, then we have
	$R_{k, z_i} \frac{1}{f} = -\frac{R_{k, z_i} f}{f^2}$. Applying the product rule \eqref{eq:productrarising} yields
	\begin{align*}
		\tilde{\partial}_{1 - k + j} \tilde{\partial}_{1 - k}^j H(Z)
		&= \sum_{i = 0}^{j + 1} \binom{j + 1}{i} \frac{R_{2 - k, z_1}^{j + 1 - i} R_{-k, z_2}^i H_1(Z)}{(R_{-1, z_1} q_Z(X))^{j + 1 - i} (-R_{-1, z_2} q_Z(X))^i} d z_1 \wedge d \mu(z_2) \\
		&- \sum_{i = 1}^j i\binom{j}{i} \frac{R_{2 - k, z_1}^{j - i} R_{-k, z_2}^i H_1(Z) (-R_{-1, z_1}R_{-1, z_2} q_Z(X))}{(R_{-1, z_1} q_Z(X))^{j + 1 - i} (-R_{-1, z_2} q_Z(X))^{i+1}} d z_1 \wedge d \mu(z_2) \\
		&+ \sum_{i = 0}^{j-1} (j - i) \binom{j}{i} \frac{R_{2 - k, z_1}^{j - i} R_{-k, z_2}^i H_1(Z) (-R_{-1, z_1}R_{-1, z_2} q_Z(X))}{(R_{-1, z_1} q_Z(X))^{j + 1 - i} (-R_{-1, z_2} q_Z(X))^{1+i}} d z_1 \wedge d \mu(z_2) \\
		&+ \tilde{\partial}_{1 + j - k} \left(G(Z) (R_{-1, z_1} R_{-1, z_2} q_Z(X)) d z_1 \wedge d \mu(z_2)\right).
	\end{align*}
	A straightforward calculation gives the result.
\end{proof}

Analogously, consider a form
\begin{align*}H(Z) &= H_2(Z) d \mu(z_1) \wedge d z_2 \\
&= \frac{1}{2} (\partial_{-1} q_Z(X)) \wedge \tilde{H}_1(Z) - \frac{1}{2} (\partial_{-1} q_Z(X)) \wedge \tilde{H}_2(Z) \in \calA^{2, 1}(\Gamma_X \bs \IH^2, \calL_{-k})
\end{align*}
with
\begin{align*}
\tilde{H}_1(Z) &= \frac{y_1^{-2} H_2(Z)}{R_{-1, z_1} q_Z(X)} d \overline{z}_1 \wedge d z_2 \in \calA^{1, 1}(\Gamma_X \bs \IH^2, \calL_{1-k}), \\
\tilde{H}_2(Z) &= \frac{H_2(Z)}{R_{-1, z_2} q_Z(X)} d \mu(z_1) \in \calA^{1, 1}(\Gamma_X \bs \IH^2, \calL_{1-k})
\end{align*}
and define $\tilde{\partial}_{1 - k} H(Z) = \partial_{1 - k}(\tilde{H}_1(Z) + \tilde{H}_2(Z))$. Then we have
\begin{align*}
	\tilde{\partial}_{1-k}^j H(Z)
	&= \sum_{i = 0}^j \binom{j}{i} \frac{R_{-k, z_1}^{j - i} R_{2 - k, z_2}^i H_2(Z)}{(R_{-1, z_1} q_Z(X))^{j - i} (-R_{-1, z_2} q_Z(X))^i} d \mu(z_1) \wedge d z_2\\
	&\quad + G(Z) (R_{-1, z_1} R_{-1, z_2} q_Z(X)) d \mu(z_1) \wedge d z_2 \in \calA^{2, 1}(\Gamma_X \bs \IH^2, \calL_{j-k}),
\end{align*}
where $G(Z)$ has no singularities along $T_X$.

\begin{Lemma}\label{lem:DiagonalRaisingOperator}
	Let $H(Z) : \IH^2 \to \IC$ be a smooth function. On the diagonal $(z,z)$ we have
	$$R_{k_1 + k_2, z}^j H(z, z) = \sum_{i = 0}^{j} \binom{j}{i} (R_{k_1, z_1}^{j - i} R_{k_2, z_2}^i H)(z, z).$$
\end{Lemma}

\begin{proof}
	This is a direct application of the chain rule.
\end{proof}

\begin{Theorem}\label{thm:AlgebraicCycleIntegrals}
	Let $H(Z)=H_1(Z) \, d z_1 \wedge d \mu(z_2) + H_2(Z) \, d \mu(z_1) \wedge dz_2 \in \calA^{2, 1}(\Gamma \bs \IH^2, \calL_{-k})$ be real analytic in $B_\varepsilon(T_X)$ and cuspidal. Then
	\begin{align*}
		\lim_{\varepsilon \to 0} \int_{\Gamma_X \bs \partial B_\varepsilon(T_X)}  \frac{H(Z)}{q_Z(X)^k} &= -\frac{8\pi \lvert Q(X) \rvert^{-k/2}}{4^k(k - 1)!}\int_{\Gamma_X \bs \H} R_{2-2k}^{k-1}\big( H_1|_{T_X}(z)  -H_2|_{T_X}(z)\big) d \mu(z),
	\end{align*}
	with the restriction to $T_X$ as defined in Proposition~\ref{prop:Restriction}.
\end{Theorem}

\begin{proof}
	For simplicity we assume that $H_2 = 0$. Let $\tilde{H}_1(Z) $ and $\tilde{H}_2(Z)$ be as in \eqref{eq:tildeh1hand2}. We first prove the assertion for $X = e_4$, so that $T_X$ is the diagonal.
	We have
	\[
	\partial \left(\frac{\tilde{H}_1(Z) + \tilde{H}_2(Z)}{(1 - k) q_Z(X)^{k - 1}}\right)
	= \frac{H(Z)}{q_Z(X)^k} + \frac{\tilde{\partial}_{1-k} H(Z)}{2(1 - k) q_Z(X)^{k - 1}}.
	\]
	Thus, by Stokes' Theorem and Corollary \ref{cor:localCalculations}
	\begin{align*}
		\lim_{\varepsilon \to 0} \int_{\Gamma_X \bs \partial B_\varepsilon(T_X)}  \frac{H(Z)}{q_Z(X)^k}
		&= \lim_{\varepsilon \to 0} \int_{\Gamma_X \bs \partial B_\varepsilon(T_X)}  \frac{\tilde{\partial}_{1-k} H(Z) + \overline{\partial} (\tilde{H}_1(Z) + \tilde{H}_2(Z))}{2(k - 1) q_Z(X)^{k - 1}} \\
		&= \lim_{\varepsilon \to 0} \int_{\Gamma_X \bs \partial B_\varepsilon(T_X)}  \frac{\tilde{\partial}_{1 - k} H(Z)}{2(k - 1) q_Z(X)^{k - 1}}.
	\end{align*}
	Here, the boundary terms vanish, because $H$ is cuspidal. Inductively we obtain using Corollary \ref{cor:localCalculations} and Lemma \ref{lem:DifferentialOperatorCalculations} and that $R_{-1, z_1} R_{-1, z_2} q_Z(X)$ vanishes on $T_X$
	\begin{align*}
		&\lim_{\varepsilon \to 0} \int_{\Gamma_X \bs \partial B_\varepsilon(T_X)}  \frac{H(Z)}{q_Z(X)^k} \\
		&= \lim_{\varepsilon \to 0} \int_{\Gamma_X \bs \partial B_\varepsilon(T_X)}  \frac{\tilde{\partial}_{1 - k}^{k - 1} H(Z)}{2^{k - 1}(k - 1)! q_Z(X)} \\
		&= \frac{2 \pi i}{2^{k - 1} (k - 1)!} \int_{\Gamma_X \bs T_X} \sum_{i = 0}^{k - 1} \binom{k - 1}{i} \frac{R_{2 - k, z_1}^{k - 1 - i} R_{-k, z_2}^i H_1(Z)}{(R_{-1, z_1} q_Z(X))^{k - 1 - i} (-R_{-1, z_2} q_Z(X))^i} d \mu(z_2).
	\end{align*}
	Using $R_{-1, z_1} q_Z(e_4) = 2i = - R_{-1, z_2} q_Z(e_4)$ on the diagonal $T_X$, this is equal to
	\begin{align*}
		&\frac{2 \pi i}{(4i)^{k - 1} (k - 1)!} \int_{\Gamma_X \bs T_X} \sum_{i = 0}^{k - 1} \binom{k - 1}{i} (R_{2 - k, z_1}^{k - 1 - i} R_{-k, z_2}^i H_1)(Z) d \mu(z_2) \\
		&= -\frac{8\pi}{(4 i)^{k} (k - 1)!} \int_{\Gamma_X \bs \H} R_{2 - 2k}^{k - 1} \big( H_1(z,z)\big) d \mu(z),
	\end{align*}
	where we used Lemma \ref{lem:DiagonalRaisingOperator}. For arbitrary $X$ we can reduce to the case above by choosing some $\gamma \in \SL_2(\IR)^2$ with  $\gamma^{-1}.X = \sqrt{|Q(X)|} \, e_4$.
	This shows the result for $H_2 = 0$. The different sign for the case with $H_2$ comes from Corollary \ref{cor:localCalculations}.
\end{proof}

\subsection{Regularized inner products}
\label{subsec:RegularizedInnerProducts}

We can now define the regularized inner product in \eqref{eq:preliminaryInnerProduct} and evaluate it for functions $G$ and $H$ with certain types of singularities along algebraic and real cycles.

Let $G \in \calL_k(\Gamma \bs \IH^2)$ be meromorphic with singularities of type
$$\sum_{\substack{X \in L' \\ Q(X) < 0}} c(G, X) \sum_{\gamma \in \Gamma_X \bs \Gamma} \frac{1}{q_Z(\gamma^{-1}. X)^k}$$
along some algebraic cycles $T_X$. Moreover, assume that $G$ is cuspidal. Let
$$U_\varepsilon = \bigcup_{\substack{X \in L' \\ c(G, X) \neq 0}} B_\varepsilon(T_X)$$
be an $\varepsilon$-neighbourhood of the singularities of $G$. 
\begin{Definition}\label{def:RegularizedInnerProduct}
We define the regularized inner product of a meromorphic modular form $G$ with a cusp form $F$ of weight $k$ as
\[
\int_{\Gamma \bs \IH^2}^{\reg} G(Z) \overline{F(Z)} (y_1 y_2)^k d\mu(Z) = \lim_{\varepsilon \to 0} \int_{\Gamma \bs (\IH^2 \setminus U_\varepsilon)} G(Z) \overline{F(Z)} (y_1 y_2)^k d\mu(Z).
\]
\end{Definition}

Assume additionally that $F(Z) = \xi_{-k,Z} H(Z)$ for a cuspidal real analytic $(2, 1)$-form $H \in \calA^{2, 1}(\Gamma \bs \IH^2, \calL_{-k})$ with singularities of type
\begin{align}
\begin{split}\label{eq:RealAnalyticSingularities}
\sum_{\substack{Y \in L' \\ Q(Y) > 0}} c(H, Y) \sum_{\gamma \in \Gamma_Y \bs \Gamma}&\left( y_1^{k-2}y_2^{k}\frac{\overline{q_Z(\gamma^{-1}.Y)}^{1-k}}{\overline{p_Z(\gamma^{-1}. Y)}} \, d z_1 \wedge d \mu(z_2) \right.  \\
&\qquad \left.+ y_1^{k}y_2^{k-2}\frac{\overline{q_Z(\gamma^{-1}.Y)}^{1-k}}{\frac{y_1}{y_2}p_{Z}(\gamma^{-1}. Y)} \, d \mu(z_1) \wedge d z_2\right)
\end{split}
\end{align}
along certain real analytic cycles $C_Y$. As before, we write
\[
H(Z) = H_1(Z) \ d z_1 \wedge d \mu(z_2) + H_2(Z)  \ d \mu(z_1) \wedge d z_2.
\]

\begin{Theorem} \label{thm:innerProductEvaluation}
Let $G$ and $H$ be as above and assume that the singularities of $G$ and the singularities of $H$ do not intersect. Then we have
\[
\int_{\Gamma \bs \IH^2}^{\reg} G(Z) \overline{\xi_{-k, Z} H(Z)} (y_1 y_2)^{k} d\mu(Z)
= \delta_{C(H)}(G) + \delta_{T(G)}(H),
\]
where 
\begin{align*}
\delta_{C(H)}(G) = -\frac{16\pi i}{4^k} \sum_{\substack{Y \in L' \\ Q(Y) > 0}} \frac{c(H, Y)}{Q(Y)^{k-1}} \int_{\Gamma_Y \bs C_Y} G(Z) q_Z(Y)^{k - 2} d z_1 \wedge d z_2,
\end{align*}
and
\begin{align*}
\delta_{T(G)}(H) =&-\frac{8\pi }{4^k(k - 1)!} \sum_{\substack{X \in L' \\ Q(X) < 0}} \frac{c(G, X)}{\lvert Q(X) \rvert^{k/2}}\int_{\Gamma_X \bs \H} R_{2-2k}^{k-1} \big(H_1|_{T_X}(z)-H_2|_{T_X}(z)\big) d \mu(z),
\end{align*}
with the restriction to $T_X$ as defined in Proposition~\ref{prop:Restriction}.
\end{Theorem}

\begin{proof}
Using Stokes' Theorem, the cuspidality of $G$ and $H$ and the meromorphicity of $G$, we can rewrite the integral to
\begin{align*}
 &\sum_{\substack{X \in L' \\ Q(X) < 0}} c(G, X) \lim_{\varepsilon \to 0} \int_{\Gamma_X \bs \partial B_{\varepsilon}(T_X)} \frac{H(Z)}{q_Z(X)^k} +\sum_{\substack{Y \in L' \\ Q(Y) > 0}} c(H, Y) \\
 & \times \lim_{\varepsilon \to 0} \int_{\Gamma_Y \bs \partial B_{\varepsilon}(C_Y)} G(Z) \left(y_1^{k-2}y_2^{k}\frac{\overline{q_Z(Y)}^{1-k}}{\overline{p_Z(Y)}} d z_1 \wedge d \mu(z_2)+y_1^{k}y_2^{k-2}\frac{\overline{q_Z(Y)}^{1-k}}{\frac{y_1}{y_2}p_{Z}(Y)} d \mu(z_1) \wedge d z_2\right).
\end{align*}
Lemma \ref{lem:RealCycleIntegrals} and Theorem \ref{thm:AlgebraicCycleIntegrals} yield the result.
\end{proof}

Finally, we show that the meromorphic modular forms $\omega_X^{\mero}$ are orthogonal to cusp forms with respect to the regularized Petersson inner product.

\begin{Theorem} \label{thm:orthogonalToCuspForms}
For every real-analytic cusp form $F : \IH^2 \to \IC$ we have
$$\int_{\Gamma \bs \IH^2}^{\reg} F(Z) \overline{\omega_X^{\mero}(Z)} (y_1 y_2)^k d\mu(Z) = 0.$$
\end{Theorem}

\begin{proof}
By \eqref{eq:OmegaXiPreimage} we have $\xi_{-k, Z} \Omega_X^{\mero} = -2(k-1) \omega_X^{\mero}$. Using Lemma \ref{lem:StokesForXi} and the usual unfolding method we have 
\begin{align*}
&2(k-1) \int_{\Gamma \bs \IH^2}^{\reg} F(Z) \overline{\omega_X^{\mero}(Z)} (y_1 y_2)^k d\mu(Z) \\
&= \lim_{\varepsilon \to 0} \int_{\Gamma_X \bs \partial B_{\varepsilon}(T_X)} F(Z) \left(y_1^{k - 2} y_2^k\frac{\overline{q_Z(X)}^{1 - k}}{\overline{p_Z(X)}}  d z_1 \wedge d \mu(z_2)
+ y_1^{k} y_2^{k-2} \frac{\overline{q_Z(X)}^{1 - k}}{\frac{y_1}{y_2} p_Z(X)}  d \mu(z_1) \wedge d z_2\right).
\end{align*}
Now apply  identity\eqref{eq:VanishingIntegral} to the complex conjugation of the integrals.
\end{proof}

\subsection{A current equation and the relation to the Oda lift}
\label{subsec:CurrentEquation}

Let $H$ be a real analytic $(2, 1)$-form with values in the line bundle of modular forms of weight $-k$ and singularities of type \eqref{eq:RealAnalyticSingularities}. Then we can consider $H$ as a distribution on rapidly decreasing forms in $\calA^{2,1}(\Gamma \bs \IH^2, \calL_{-k})$ via
$$[H](G) = \int_{\Gamma \bs \IH^2} G \wedge \overline{*}_{-k} H.$$
We define the $\xi_{-k, Z}$ operator on these distributions via
$$(\xi_{-k, Z}[H])(G) = - [H](\xi_{k, Z} G),$$
where $G \in \calA^0(\Gamma \bs \IH^2, \calL_{k})$ is rapidly decreasing.

\begin{Theorem}\label{thm:DistributionEquation}
	With notation as above we have
$$\xi_{-k, Z} [H] = [\xi_{-k, Z} H] - \delta_{C(H)},$$
where $\delta_{C(H)}$ is given as in Theorem~\ref{thm:innerProductEvaluation}.
\end{Theorem}

\begin{proof}
We have
\[
\xi_{-k, Z}[H](G)
	= -\int_{\Gamma \bs \IH^2} \xi_{k, Z} G(Z) \wedge \overline{*}_k H(Z).
\]
By Lemma \ref{lem:StokesForXi} this equals
\[
\int_{\Gamma \bs \IH^2} G(Z) \overline{\xi_{-k, Z} H(Z)} (y_1 y_2)^k d \mu(Z) - \int_{\Gamma \bs \IH^2} d (G(Z) H(Z)).
	\]
The first integral is given by $[\xi_{-k,Z} H](G)$. We turn to the evaluation of the second integral. Using Stokes' Theorem and the rapid decay of $G$ we obtain
\begin{align*}
&\int_{\Gamma \bs \IH^2} d (G(Z) H(Z))
= \sum_{\substack{Y \in L' \\ Q(Y) > 0}} c(h, Y) \\
&\times  \lim_{\varepsilon \to 0} \int_{\Gamma_Y \bs \partial B_{\varepsilon}(C_Y)} G(Z)\left( y_1^{k - 2} y_2^k \frac{\overline{q_Z(X)}^{1 - k}}{\overline{p_Z(X)}} d z_1 \wedge d \mu(z_2)
+ y_1^k y_2^{k - 2} \frac{\overline{q_Z(X)}^{1 - k}}{\frac{y_1}{y_2} p_Z(X)}  d \mu(z_1) \wedge d z_2\right).
\end{align*}
The application of Lemma \ref{lem:RealCycleIntegrals} gives the result.
\end{proof}

\begin{Remark}
	As a special case of Theorem \ref{thm:DistributionEquation} we have
	$$\langle \omega_{Y_1}^{\cusp}, \omega_{Y_2}^{\cusp} \rangle = \frac{2 \pi i}{k - 1} \int_{\Gamma_{Y_2} \bs C_{Y_2}} \omega_{Y_1}^{\cusp}(z, -\overline{z}) q_Z(Y_2)^{k - 2} dz_1 \wedge dz_2.$$
	Integrals of this form have been calculated in some cases by \cite{ZagierCycleIntegralsOmega} and \cite{ZagierTracesHecke} and are given by finite sums involving Hurwitz class numbers. In \cite{katok}, Katok calculates similar integrals in signature $(2, 1)$ and derives two rational structures on the spaces of cusp forms. 
\end{Remark}

Using Theorem \ref{thm:ThetaLiftAndXiOperator} we obtain the following current equation for the Millson and the Doi-Naganuma lift, which is the analogue of \cite[Theorem 8.4]{crawfordfunke}.

\begin{Corollary}
For $g \in H^+_{2-k, L^-}$ we have
$$\xi_{-k, Z} [\Phi_L^{\Millson}(g)] = [\Phi_L^{\DN}(\xi_{2-k, \tau} g)] - \delta_{C(\Phi_L^{\Millson}(g))}.$$
\end{Corollary}

Following \cite{oda} we define the Oda lift of a cusp form $F \in \calL_k(\Gamma \bs \IH^2)$ by
$$\Phi_L^{\Oda}(F,\tau) = \int_{\Gamma \bs \IH^2} F(Z) \Theta_L^{\DN}(\tau, Z) (y_1 y_2)^k d\mu(Z).$$
Then we have the following analogue of \cite[Lemma 8.7]{crawfordfunke}.

\begin{Corollary}
For a harmonic Maass form $g \in H^+_{2 - k, L^-}$ and a cusp form $F \in \calL_k(\Gamma \bs \IH^2)$ we have
\begin{align*}
\int_{\Mp_2(\IZ) \bs \IH} \left\langle \Phi_L^{\Oda}(F,\tau), \xi_{2 - k, \tau} g(\tau) \right\rangle v^k d \mu(\tau) = \frac{1}{4i} \delta_{C(\Phi_L^{\Millson}(g))}(F).
\end{align*}
\end{Corollary}

Choosing $g$ as a Maass-Poincar\'e series defined in Section~\ref{subsec:MaassPoincareSeries} yields the Fourier expansion of the Oda lift.

\section{Cycle integrals of meromorphic Hilbert modular forms}
\label{sec:CycleIntegrals}

In this section, we prove an explicit formula for the traces of cycle integrals of the meromorphic modular forms $\omegamero{X}$ defined in Section~\ref{subsec:MeromorphicModularForms}. Additionally, we show that, up to taking linear combinations of these cycle integrals, these traces are rational.

Recall that, for $k \geq 4$ even and $X \in V$ with $Q(X) < 0$, we defined a meromorphic modular form of weight $k$ for $\Gamma$ by
\[
\omegamero{X}(Z) = \sum_{\gamma \in \Gamma_X \backslash \Gamma} q_Z(\gamma^{-1} . X)^{-k}
\]
and noted that it vanishes at the cusps and has poles at the algebraic cycle $T_X$ and its $\Gamma$-translates. For $\beta \in L'/L$ and $m > 0$ we define its $m$-th trace of cycle integrals by
\begin{align*}
	\tr_{m,\beta}(\omega_X^{\mero}) = \sum_{Y \in \Gamma \backslash L_{m,\beta}}\int_{\Gamma_Y \backslash C_Y} \omegamero{X}(Z)q_Z(Y)^{k-2} dz_1 \wedge dz_2.
\end{align*}

For simplicity, we will assume that the cycles $T_X$ and $C_Y$ for $Y \in L_{m,\beta}$ do not intersect. By Lemma \ref{lemma intersection cycles}, this means that we have $(X,Y) \neq 0$ for all $Y \in L_{m,\beta}$. 

As before, we define the sublattices
\[
P = L \cap (\Q X)^\perp, \qquad N = L \cap (\Q X),
\]
which have signature $(2,1)$ and $(0,1)$, respectively. The traces of cycle integrals of $\omega_X^{\mero}$ can be computed explicitly as detailed in the following theorem.

\begin{Theorem} \label{thm:ExplicitFormula}
	If the lattice $P = L \cap X^\perp$ is anisotropic, then we have $\tr_{m,\beta}(\omega_X^{\mero}) = 0$. Otherwise, we have
	\begin{align*}
		\tr_{m,\beta}(\omega_X^{\mero}) &= \frac{\sqrt{2}\pi i}{C}\sum_{\ell \in \Gamma_X \backslash \Iso(P)}\alpha_\ell \left(\CT\left\langle f_{m, \beta,P \oplus N}(\tau), \overline{\left[\Theta_{K_\ell,\frac{3}{2}}(\tau), \widetilde{\Theta}_{N^-}^{*,+}(\tau) \right]}_{\frac{k-2}{2}}\right\rangle \right. \\
		&\qquad \qquad \qquad \left. -\int_\calF^\text{reg} \left\langle \left[\Theta_{K_\ell,\frac{3}{2}}(\tau), \widetilde{\Theta}_{N^-}^*(\tau) \right]_{\frac{k-2}{2}}, \xi_{2-k}f_{m, \beta,P \oplus N}(\tau) \right\rangle v^k d\mu(\tau)\right)
	\end{align*}
	with the rational constant $C = \binom{k-3}{k/2-1}(k-1)16 Q(X)^{k/2}$. Here, $\CT$ denotes the constant term of a holomorphic $q$-series, $\langle \cdot,\cdot \rangle$ denotes the inner product on $\C[L'/L]$, and $[\cdot,\cdot]_{n}$ is the Rankin-Cohen bracket in weights $\frac{3}{2}$ and $\frac{1}{2}$ as defined in \eqref{eq:RankinCohen}.
\end{Theorem}

\begin{Remark}
We also remind the reader of the following notation we used in the theorem: $\Iso(P)$ is the set of isotropic lines in $P \otimes \Q$, $\alpha_\ell \in \R_{> 0}$ is the width of the cusp corresponding to $\ell$, and $\Theta_{K_\ell,\frac{3}{2}}$ is the cuspidal weight $3/2$ theta series associated with $K_\ell$ (see Section~\ref{subsec:ThetaFunctions21}). Finally, $f_{m,\beta,P\oplus N}$ is the image under the map $f_{P \oplus N}$ in \eqref{eq:ArrowDown} of the Maass-Poincar\'e series $f_{m,\beta}$ of weight $2-k$ for $\rho_{L^-}$ defined in Section~\ref{subsec:MaassPoincareSeries}, $\widetilde{\Theta}_{N^-}^{*}$ is a harmonic Maass form of weight $1/2$ for $\rho_N$ which maps to the weight $3/2$ unary theta function $\Theta_{N^-}^*$ (defined in \eqref{eq:UnaryThetaFunctions}) under $\xi_{1/2}$, and the superscript $+$ denotes its holomorphic part. 
\end{Remark}

\begin{proof}
	We start by considering the regularized inner product
	\begin{align*}
		\left\langle \omegamero{X},\omegacusp{Y}\right\rangle_{\Pet}^{\reg} = \int_{\Gamma \backslash \H^2}^{\reg} \omegamero{X}(Z) \overline{\avgomegacusp{\beta}{m}(Z)}(y_1y_2)^k d\mu(Z).
	\end{align*}
	On the one hand this vanishes by Theorem~\ref{thm:orthogonalToCuspForms}. On the other hand, since $\xi_{-k,Z}\Omega_{m,\beta}^{\cusp} = -2(k-1)\omega_{m,\beta}^{\cusp}$, it it is a constant multiple of
	\[\delta_{C(\Omega_{m,\beta}^{\cusp})}(\omegamero{X}) + \delta_{T(\omegamero{X})}(\avgOmegacusp{m}{\beta})\]
	by Theorem~\ref{thm:innerProductEvaluation}. This implies
	\[
	\tr_{m,\beta}(\omegamero{X}) = \frac{i}{2(k-1)!}\frac{m^{k-1}}{Q(X)^{k/2}}\Omega_{m,\beta}^{\cusp}(T_X),
	\]
	where $\Omega_{m,\beta}^{\cusp}(T_X)$ denotes the evaluation of $\Omega_{m,\beta}^{\cusp}$ at $T_X$ as in Definition~\ref{def:EvaluationAtTX}. To compute this, we use the fact that $\avgOmegacusp{m}{\beta}$ can be written as the image of the Poincaré series $f_{m,\beta}(\tau)$ via the Millson theta lift $\Phi_L^{\Millson}$ by Theorem \ref{thm:unfoldpc}, that is
	\begin{align*}
		&\avgOmegacusp{m}{\beta}(Z) = \frac{\pi}{2(4m)^{k-1}}\Phi_{L}^{\Millson}(f_{m,\beta},Z).
	\end{align*}
	After plugging in the definition of the lift, exchanging the order of integration, and moving the raising operator into the inner integral, we obtain
	\begin{align*}
		\tr_{m,\beta}(\omega_X^{\mero})&= \frac{\pi i}{4^{k}(k-1)!Q(X)^{k/2}}\int_\calF^\text{reg} \left\langle f_{m, \beta}(\tau) , \overline{\Theta_{L}^{\Millson}(\tau,T_X)} \right\rangle v^{2-k} d\mu(\tau).
	\end{align*}
	We plug in the evaluation of $\Theta_{L}^{\Millson}(\tau,Z)$ at $T_X$ from Theorem~\ref{thm:EvaluationMillsonTX} and obtain
	\begin{align*}
		\tr_{m,\beta}(\omega_X^{\mero}) &= \frac{(-1)^{k/2+1}\pi i \sqrt{2}^{k+1}}{4^{k}(k-1)Q(X)^{k/2}} \\
		&\quad \times \sum_{\ell \in \Gamma_X \backslash \Iso(P)}\alpha_\ell\int_\calF^\text{reg} \left\langle f_{m, \beta}(\tau) , \left(\overline{\Theta_{K_\ell,k-1/2}(\tau) \otimes \overline{\Theta_{N^-}^*(\tau)}} \right)^L \right\rangle v^{\frac{3}{2}} d\mu(\tau).
	\end{align*}	
	Since the operators $g^{L}$ and $f_{P \oplus N}$ are adjoint to each other, we can erase the superscript $L$ and add a subscript $P \oplus N$ to $f_{m,\beta}$. Now we compute the integral for fixed $\ell \in \Iso(P)$. By Lemma~\ref{lem:RaisingUnaryThetas} we have
	\[\Theta_{K_{\ell}, k-\frac{1}{2}}(\tau) = \left(-2\pi\right)^{\frac{2-k}{2}} R_{\frac{3}{2},\tau}^{\frac{k-2}{2}} \Theta_{K_{\ell}, \frac{3}{2}}(\tau),\]
	so the integral can be written as
\begin{align*}
		(-2\pi)^{\frac{2-k}{2}}\int_\calF^\text{reg} \left\langle f_{m, \beta,P \oplus N}(\tau) , \left(\overline{R_{\frac{3}{2},\tau}^{\frac{k-2}{2}}\Theta_{K_\ell,\frac{3}{2}}(\tau) \otimes \overline{\Theta_{N^-}^*(\tau)}} \right) \right\rangle v^{\frac{3}{2}} d\mu(\tau).
	\end{align*}
	Let $\widetilde{\Theta}^*_{N^-}(\tau)$ be a harmonic Maass form of weight of $1/2$ for $\rho_{N}$ which maps to $\Theta^*_{N^-}(\tau)$ under the $\xi$-operator, and let $\widetilde{\Theta}^{*,+}_{N^-}$ denote its holomorphic part. Using Proposition \ref{rankin cohen brackets} we rewrite the previous expression as 
	\begin{align*}
		(-2\pi)^{\frac{2-k}{2}}(4\pi)^{\frac{k-2}{2}}\frac{\Gamma(\frac{k}{2})\Gamma(\frac{1}{2})}{\Gamma(\frac{k}{2}-\frac{1}{2})}\int_\calF^\text{reg} \left\langle f_{m, \beta,P \oplus N}(\tau), \overline{L\left[\Theta_{K_\ell,\frac{3}{2}}(\tau), \widetilde{\Theta}_{N^-}^*(\tau) \right]}_{\frac{k-2}{2}} \right\rangle  d\mu(\tau),
	\end{align*}
	where we used that $\Theta_{N^-}^*(\tau)$ is holomorphic on $\H$. 
	We can simplify the constant in front to $(-1)^{k/2+1}2^{3k/2-4}/\binom{k-3}{k/2-1}$ using the Legendre duplication formula.

Using Stokes' Theorem in the form given in \cite[Lemma~2.1]{brikavia}, implies that the integral in the last expression equals
	\begin{align*}
		&- \int_{\partial \calF} \left\langle f_{m, \beta,P \oplus N}(\tau) , \overline{\left[\Theta_{K_\ell,\frac{3}{2}}(\tau), \widetilde{\Theta}_{N^-}^*(\tau) \right]}_{\frac{k-2}{2}} \right\rangle dz \\
		&-\int_\calF^\text{reg} \left\langle \overline{\xi_{2-k}f_{m, \beta,P \oplus N}(\tau)} , \overline{\left[\Theta_{K_\ell,\frac{3}{2}}(\tau), \widetilde{\Theta}_{N^-}^*(\tau) \right]}_{\frac{k-2}{2}} \right\rangle v^k d\mu(\tau).
	\end{align*}
	The first integral reduces to an integral on the horizontal boundary of the standard regularized fundamental domain since all other contributions cancel out due to the modularity of the integrand. Taking the limit as this horizontal line goes to infinity picks out the constant term of	the Fourier expansion of the expression inside the integral. Namely, we arrive to 
	\begin{align*}
		&\CT\left\langle f_{m, \beta,P \oplus N}(\tau) , \overline{\left[\Theta_{K_\ell,\frac{3}{2}}(\tau), \widetilde{\Theta}_{N^-}^{*,+}(\tau) \right]}_{\frac{k-2}{2}}\right\rangle \\
		&-\int_\calF^\text{reg} \left\langle \left[\Theta_{K_\ell,\frac{3}{2}}(\tau), \widetilde{\Theta}_{N^-}^*(\tau) \right]_{\frac{k-2}{2}}, \xi_{2-k}f_{m, \beta,P \oplus N}(\tau) \right\rangle v^k d\mu(\tau).
	\end{align*}

	Notice that only the holomorphic part of $\widetilde{\Theta}_{N^-}^{*}$ remained in the boundary integral since the non-holomorphic part decreases square-exponentially, and hence compensates for any terms appearing in the Rankin-Cohen bracket. Gathering all the constants, we obtain the stated formula.
\end{proof}

\begin{Theorem}\label{thm:MainRationalityResult}
	Let $f \in M_{2-k,L^-}^!$ be a weakly holomorphic modular form with rational Fourier coefficients $c_{f}(-m,\beta)$ for $m < 0$. Moreover, assume that the cycles $T_X$ and $C_Y$ for $Y \in L_{m,\beta}$ do not intersect if $c_f(-m,\beta)\neq 0$. Then, the linear combination
	\[\sum_{\beta \in L'/L} \sum_{m >0} c_f(-m,\beta) \tr_{m,\beta}(\omegamero{X})\]
	of traces of cycle integrals is a rational multiple of $\pi i$.
\end{Theorem}

\begin{proof}
	First notice that the Petersson inner product on the right-hand side of the formula in Theorem~\ref{thm:ExplicitFormula} vanishes since $f$ is weakly holomorphic. Hence, by writing $f= \frac{1}{2} \sum_{\beta \in L'/L} \sum_{m >0} c_f(-m,\beta)f_{m,\beta}$ as a linear combination of Maass-Poincar\'e series, we see that the linear combination of traces of cycle integral in the theorem is a rational multiple of 
	\begin{align}\label{eq:FinalFormula}
	\sqrt{2}\pi i\sum_{\ell \in \Gamma_X \backslash \Iso(P)}\alpha_\ell\bigg(\CT\left( f_{P \oplus N}(\tau) \cdot \left[\Theta_{K_\ell,\frac{3}{2}}(\tau), \widetilde{\Theta}_{N^-}^{*,+}(\tau) \right]_{\frac{k-2}{2}}\right).
	\end{align}
	By Remark \ref{rem:alphaell}, the width of the cusp $\alpha_\ell$ is a rational multiple of $\sqrt{\lvert K_\ell' / K_\ell \rvert / 2}$. In particular, the $\sqrt{2}$ cancels and we are left to show that the rest of the terms cancel the factor $\sqrt{\lvert K_\ell' / K_\ell \rvert}$.
	Therefore, we need to investigate the algebraic properties of the Fourier coefficients of the harmonic Maass forms appearing in this formula.
	
	Since the space $M_{2-k,L^-}^!$ of weakly holomorphic modular forms has a basis of forms with rational Fourier coefficients by a result of McGraw \cite{mcgraw}, the assumption that the coefficients $c_f(-m,\beta)$ for $m > 0$ are rational implies that \emph{all} coefficients of $f$ are rational. 
	
	Let us now consider $\widetilde{\Theta}_{N^-}^{*}(\tau)$. Inspecting the definition \eqref{eq:UnaryThetaFunctions}, we see that the unary theta function $\Theta_{N^-}^*(\tau)$ is a rational multiple of $\sqrt{|Q(X)|}\theta_{\mathcal{N}}(\tau,1)$, where $\theta_{\mathcal{N}}(\tau;1)$ is the unary theta function defined in the introduction of \cite{lischwagenscheidt} and $\mathcal{N} \in \N$ is such that $N^-$ is isomorphic to $(\Z,\mathcal{N}x^2)$. Since the lattice $N^-$ is generated by a rational multiple of $X$, we see that $|Q(X)|$ is a rational square multiple of $\mathcal{N}$. Hence, by \cite[Theorem~1.1]{lischwagenscheidt}, we can choose the harmonic Maass form $\widetilde{\Theta}_{N^-}^{*}(\tau)$ such that its holomorphic part has rational coefficients.

	It remains to investigate the algebraic nature of the Fourier coefficients of the unary theta series
	\[
	\Theta_{K_\ell,\frac{3}{2}}(\tau) = \sum_{\lambda \in K_\ell'}(\lambda,w_\ell)e(Q(\lambda)\tau)\e_{\lambda + K_\ell}.
	\]
	We show that its coefficients are rational multiples of $1 / \sqrt{\lvert K_\ell' / K_\ell \rvert}$.
	Recall that $w_\ell$ is a generator of the real quadratic space $W_\ell(\R) = K_\ell \otimes \R$ with $(w_\ell,w_\ell) = 1$. We can write it as $w_\ell = \kappa_\ell/|\kappa_\ell|$, where $\kappa_\ell \in K_\ell$ is primitive, and $|x| = \sqrt{(x,x)}$. Notice that $(\lambda,\kappa_\ell) \in \Z$ for $\lambda \in K_\ell'$, so that the coefficients are rational multiples of $1/|\kappa_\ell| = 1/\sqrt{(\kappa_\ell, \kappa_\ell)} = 1/\sqrt{\lvert K_\ell' / K_\ell \rvert}$.

	Finally, note that the Rankin-Cohen bracket preserves the rationality of the Fourier coefficients. This finishes the proof.
\end{proof}

\begin{proof}[Proof of Theorem~\ref{thm:MainResultIntro}]
	We briefly explain how Theorem~\ref{thm:MainResultIntro} in the introduction follows from Theorem~\ref{thm:MainRationalityResult}. Let $L$ be the lattice coming from the ring of integers  in a number field $F = \Q(\sqrt{D})$ with discriminant $D > 0$ as in the introduction. If $D$ is an odd prime, then by \cite{bruinierbundschuh} there is an isomorphism $S_{k,L} \cong S_{k}^+(\Gamma_0(D),\chi_D)$. In particular, if the latter space is trivial, then $S_{k,L}$ is trivial, and hence every Maass-Poincar\'e series $f_{m,\beta} \in H_{2-k,L^-}^+$ is weakly holomorphic. Thus, we do not need to take linear combinations of traces of cycle integrals in Theorem~\ref{thm:MainRationalityResult} in this case. We obtain the rationality statement of Theorem~\ref{thm:MainResultIntro}.
	\end{proof}

\bibliography{bib}{}
\bibliographystyle{alpha}

\end{document}